\documentclass[reqno]{amsart}
\usepackage{amsmath,amssymb,amscd}

\providecommand{\ch}{\mathop{\rm char}\nolimits}

\DeclareMathOperator{\Spec}{Spec} \DeclareMathOperator{\id}{id}
\DeclareMathOperator{\GK}{GKdim} \DeclareMathOperator{\lex}{lex}
 \DeclareMathOperator{\ad}{ad}
\DeclareMathOperator{\Irr}{Irr} 
 \DeclareMathOperator{\hgt}{ht}
 \DeclareMathOperator{\Aut}{Aut}
\DeclareMathOperator{\PZ}{PZ} \DeclareMathOperator{\ham}{ham}
 \DeclareMathOperator{\gr}{gr}

\newcommand{\ov}{\overline}
\newcommand{\Z}{{\mathbb Z}}

\newcommand{\N}{{\mathbb N}}
\newcommand{\F}{{\mathbb F}}

\newcommand{\K}{{\mathbb K}}

\theoremstyle{plain}
\newtheorem{theorem}{Theorem}[section]
\newtheorem{lemma}[theorem]{Lemma}
\newtheorem{cor}[theorem]{Corollary}
\newtheorem{prop}[theorem]{Proposition}

\theoremstyle{definition}
\newtheorem{defn}[theorem]{Definition}
\newtheorem{defns}[theorem]{Definitions}
\newtheorem{rmk}[theorem]{Remark}

\newtheorem{example}[theorem]{Example}
\newtheorem{examples}[theorem]{Examples}

\newtheorem{notn}{Notation}

\numberwithin{equation}{section}

\author[D. A. Jordan]{David A. Jordan}

\address{Department of Pure Mathematics\\
University of  Sheffield\\
Hicks Building\\
Sheffield S3~7RH\\
UK.}

\email{d.a.jordan@sheffield.ac.uk}
\author[N. Sasom]{Nongkhran Sasom}

\address{Department of Pure Mathematics\\
University of  Sheffield\\
Hicks Building\\
Sheffield S3~7RH\\
UK.} \curraddr{
   Department of Mathematics,  Statistics and Computer Science\\
  Faculty of Science\\
 Ubon Ratchathani University\\
  Warin Chamrap\\
  Ubon Ratchathani\\
  34190\\
  Thailand}

   \email{sasomn@sci.ubu.ac.th}

  \thanks{The second
author was supported by a Royal Thai Government Scholarship.}

\title[Reversible skew Laurent polynomial rings]
{Reversible skew Laurent polynomial rings and deformations of Poisson
automorphisms}

\begin{document}
\begin{abstract}
A skew Laurent polynomial ring $S=R[x^{\pm 1};\alpha]$ is reversible
if it has a reversing automorphism, that is, an automorphism
$\theta$ of period $2$ that transposes $x$ and $x^{-1}$ and
restricts to an automorphism $\gamma$ of $R$ with
$\gamma=\gamma^{-1}$. We study invariants for reversing
automorphisms and apply our methods to determine the rings of
invariants of reversing automorphisms of the two most familiar
examples of simple skew Laurent polynomial rings, namely a
localization of the enveloping algebra of the two-dimensional
non-abelian solvable Lie algebra and the coordinate ring of the
quantum torus, both of which are deformations of Poisson algebras
over the base field $\F$. Their reversing automorphisms are
deformations of Poisson automorphisms of those Poisson algebras. In
each case, the ring of invariants of the Poisson automorphism is the
coordinate ring $B$ of a surface in $\F^3$ and the ring of
invariants $S^\theta$ of the reversing automorphism is a deformation
of $B$ and is a factor of a deformation of  $\F[x_1,x_2,x_3]$ for  a
Poisson bracket determined by the appropriate surface.
\end{abstract}

\maketitle

\section{Introduction}
\begin{notn}
Throughout  $\F$ denotes a field, $\Aut(R)$ denotes the group of
automorphisms of a ring $R$ and if $R$ is an $\F$-algebra then
$\Aut_\F(R)$ is the group of $\F$-automorphisms of $R$. Whenever we
discuss Poisson $\F$-algebras, we shall assume that $\ch \F=0$. We
denote by $\N_0$ the set of non-negative integers.
\end{notn}

If $R$ is any ring then there is an automorphism $\theta$ of the
Laurent polynomial ring $R[x^{\pm 1}]$ such that $\theta(x)=x^{-1}$
and $\theta(r)=r$ for all $r\in R$. A skew Laurent polynomial ring
$S=R[x^{\pm 1};\alpha]$, where $\alpha$ is an automorphism of $R$,
has no such automorphism unless $\alpha^2=\id_R$. However there may
exist $\theta\in \Aut(S)$, of order $2$, such that
$\theta(x)=x^{-1}$ and the restriction $\theta|_R$ is an
automorphism $\gamma$ of $R$, necessarily such that
$\gamma^2=\id_R$. We shall see, in Proposition~\ref{whenreversible},
that such an automorphism $\theta$ exists if and only if
$\gamma\alpha\gamma^{-1}=\alpha^{-1}$, in which case we say that
 $\alpha$ is $\gamma$-reversible, that $\theta$ is a {\it
reversing} automorphism and that $S$ is a {\it reversible} skew
Laurent polynomial ring. The concept of reversibility arises in
dynamical systems and the theory of flows, for example see
 \cite{devaney,devaneybook,jarczyk,henon,ofarrell}.

One of the two best known examples of simple skew Laurent polynomial
rings is the localization $V(\mathfrak{g})=\F[x^{\pm 1},y:xy-yx=x]$
at the powers of the normal element $x$ of the enveloping algebra
$U(\mathfrak g)=\F[x,y:xy-yx=x]$ of the two-dimensional non-abelian
solvable Lie algebra $\mathfrak{g}$. This is $\F[y][x^{\pm
1};\alpha]$, where,$\alpha(y)=y+1$ and it is simple provided $\ch
\F=0$. The second is the coordinate ring $W_q=\F[x^{\pm
1},y^{\pm1}:xy=qyx]=\F[y^{\pm1}][x^{\pm 1};\alpha]$ of the quantum
torus, where $q\in \F\backslash\{0\}$ and  $\alpha(y)=qy$, and is
simple provided $q$ is not  a root of unity. Both these examples are
reversible, the appropriate automorphisms $\gamma$ being such that
$\gamma(y)=-y$, for $V(\mathfrak{g})$, and $\gamma(y)=y^{-1}$,  for
$W_q$. A common approach will be used to compute the invariants for
the reversing automorphisms of $V(\mathfrak{g})$ and $W_q$,
together, in each case, with those for a reversing automorphism of
an associated reversible skew Laurent polynomial ring of which it is
a factor. For $V(\mathfrak{g})$, this is the localized homogenized
enveloping algebra $V_t(\mathfrak{g})=\F[x^{\pm 1},y,t:xy-yx=x,
xt=tx, yt=ty]$ of $\mathfrak{g}$ and for $W_q$ it is $W_Q=\F[x^{\pm
1},y^{\pm1}, Q^{\pm 1}:xy=Qyx, xQ=Qx, yQ=Qy]$, the coordinate ring
of the generic quantum torus.

Suppose now that $\ch \F=0$. If $T$ is an $\F$-algebra with a
central non-unit non-zero-divisor $t$ such that $B:=T/tT$ is
commutative then there is a Poisson bracket $\{-,-\}$ on $B$ such
that $\{\ov u,\ov v\}=\ov{t^{-1}[u,v]}$ for all $\ov u,\ov v\in B$.
In this situation, we shall follow \cite[Chapter III.5]{BGl} in
referring to $T$ as a quantization of the Poisson algebra $B$ and we
shall refer to an $\F$-algebra of the form
$T_\lambda=T/(t-\lambda)T$, where $\lambda\in \F$ is such that the
central element $t-\lambda$ is a non-unit in $T$, as a {\it
deformation} of $B$. In this sense, $V_t(\mathfrak{g})$ and
$V(\mathfrak{g})$ are, respectively, a quantization and deformation
of $\F[x^{\pm 1},y]$, with $\{x,y\}=x$, while, taking $t=Q-1$ and
$\lambda=q-1$, $W_Q$ and $W_q$ are, respectively, a quantization and
deformation of $\F[x^{\pm 1},y^{\pm 1}]$, with $\{x,y\}=xy$.

With $T, t, B$ and $\lambda$ as above, let $\theta\in \Aut_\F (T)$
be such that $\theta(t)=t$. Such an automorphism induces, in obvious
ways, a Poisson $\F$-automorphism $\pi$ of $B$ and an automorphism
$\theta_\lambda$ of $T_\lambda$. We shall refer to $\theta$ and
$\theta_\lambda$, respectively, as a quantization and a deformation
of $\pi$. The reversing automorphisms $\theta$ of
$V_t(\mathfrak{g})$ and $\theta_1$ of $V(\mathfrak{g})$ are,
respectively, a quantization and deformation of the Poisson
automorphism $\pi$ of $\F[x^{\pm 1},y]$ such that $\pi(x)=x^{-1}$
and $\pi(y)=-y$.  The ring of invariants of $\pi$ is a Poisson
subalgebra of $\F[x^{\pm 1},y]$ and is readily seen to be isomorphic
to the coordinate ring of the surface $x_1(x_2^2-4)-x_3^2=0$ in
$\F^3$.

The reversing automorphism  $\theta$ of $W_Q$ is a quantization of
the Poisson automorphism $\pi$ of the coordinate ring $\F[x^{\pm
1},y^{\pm 1}]$ of the torus such that $\pi(x)=x^{-1}$ and
$\pi(y)=y^{-1}$. We shall write $\theta_q$, rather than
$\theta_{q-1}$, for the corresponding automorphism of $W_q$ which
deforms $\pi$. The ring of invariants of $\pi$ is a Poisson
subalgebra of $\F[x^{\pm 1},y^{\pm 1}]$ and is known to be
isomorphic to the coordinate ring of the surface
$x_1x_2x_3=x_1^2+x_2^2+x_3^2-4$ in $\F^3$. For example, see
\cite[Example 3.5]{lor}, although there the base ring is $\Z$ rather
than a field.

For each of our main examples, the situation is represented in
Figure~\ref{quantPoiss} where the top row consists of Poisson
algebras and Poisson homomorphisms. Here $A:=\F[x_1,x_2,x_3]$ ,
$f\in A$ is irreducible and $B=\F[x^{\pm 1},y]$ or $\F[x^{\pm
1},y^{\pm 1}]$. The second and third rows are a quantization and a
deformation of the first. Here $g$ and $p$ are central elements of
the quantization $T$ and deformation $T_\lambda$ of $A$. In the last
two columns, $S$ and $S_\lambda$ are reversible skew Laurent
polynomial rings, either $V(\mathfrak{g})$ and ${U}(\mathfrak{g})$
or $W_Q$ and $W_q$, and $\theta$ and $\theta_\lambda$ are reversing
automorphisms. Each $j_i$ is inclusion and each $p_i, q_i$ or $d_i$
is a natural surjection. The maps $p_2$ and $p_3$ may be regarded as
a quantization and deformation of the embedding of the appropriate
surface in $\F^3$.
\begin{figure}[h]
\label{quantPoiss}
\[
\begin{CD}
 A @>p_1>> A/fA @>\simeq >> B^\pi @>j_1>> B @>\pi >> B \\
 @Aq_1AA @ Aq_2AA@Aq_3AA@Aq_4AA@Aq_4AA\\
 T @>p_2>> T/gT@>\simeq>> S^{\theta}
  @>j_2>> S@>\theta >> S \\
 @Vd_1VV @ Vd_2VV@Vd_3VV@Vd_4VV@Vd_4VV\\
 T_\lambda @>p_3>> T_\lambda/pT_\lambda@>\simeq >> S_\lambda^{\theta_\lambda} @>j_3>> S_\lambda@>\theta_\lambda >> S_\lambda\\
\end{CD}
\]
\caption{Quantization, deformation and invariants}
\end{figure}

In the case of the localized enveloping algebra, the deformation
$T_1$ is an iterated skew polynomial ring in three indeterminates
but, for the quantum torus, no such structure is apparent for the
deformation which  has been of interest elsewhere in the literature.
It arises as the cyclically $q$-deformed enveloping algebra
$U^\prime_q(so_3)$ \cite{mp1,mp2} and as the algebra determined by a
special case of the Askey-Wilson relations
\cite{terwilliger,terwilliger+}. Both deformations are examples of
algebras determined by noncommutative potentials, as described in
\cite{Ginzburg}.

The invariants for the reversing automorphisms of the localized
enveloping algebra and its homogenization are computed in
Section~\ref{invlocenv} and those for the quantum torus and generic
quantum torus are computed in Section~\ref{invqtorus}.
Section~\ref{revbasics} contains the definitions and the main
examples of reversing automorphisms together with some general
results on generators and relations for their rings of invariants.
Basic material on Poisson structures, quantization and deformation
appears in Section~\ref{PQD} while Section~\ref{filtrations}
presents some technical material, on filtrations and the Diamond
Lemma, that is applied in Sections~\ref{invlocenv} and
\ref{invqtorus}.

Some of the results of the paper appeared in the PhD thesis of the
second author \cite{nongthesis}. The study will be continued in two
papers by the first author \cite{pbps,rev2}. The notion of reversing
automorphism will be extended to other algebras, including $U(sl_2)$
and $U_q(sl_2)$, and a connection between the reversing
automorphisms of $U_q(sl_2)$ and $U^\prime_q(so_3)$
 will be exploited to determine the prime spectrum $\Spec
(U^\prime_q(so_3))$. The Poisson spectrum for certain Poisson
brackets on $A$, including the two main examples of this paper, will
be identified and, for each of those examples, it will be shown that
there is a homeomorphism between the Poisson spectrum and the
completely prime subspace of the spectrum of the corresponding
deformation $S_\lambda$.

\section{Reversing automorphisms and invariants}
\label{revbasics} \begin{defn} Let $R$ be a ring and let
$\alpha,\gamma\in \Aut(R)$ be such that $\gamma^2=\id_R$. We say
that $\alpha$ is $\gamma$\emph{-reversible} if
$\gamma\alpha\gamma^{-1}=\alpha^{-1}$. In other words $\alpha$ and
$\gamma$ provide a representation of the infinite dihedral group in
$\Aut(R)$.
\end{defn}
 It is easy to check that
$\gamma$-reversibility of $\alpha$ is equivalent to each of the
following four statements: (a) $(\alpha\gamma)^2=\id_R$;
    (b) $(\gamma\alpha)^2=\id_R$;
    (c) $\alpha=\gamma\tau$ for some $\tau\in \Aut(R)$ such that
    $\tau^2=\id_R$; (d)
     $\alpha=\tau^\prime\gamma$ for some $\tau^\prime\in \Aut(R)$ such that
    $\tau'^2=\id_R$.
\begin{prop} \label{whenreversible} Let $R$ be a ring and let $\alpha,\gamma\in\Aut(R)$ be
such that $\gamma^2=\id_R$. Let $S=R[x^{\pm1};\alpha]$. There exists
$\theta\in\Aut(S)$ such that $\theta|_R=\gamma$ and
$\theta(x)=x^{-1}$ if and only if $\alpha$ is $\gamma$-reversible.
\end{prop}
\begin{proof} Suppose that such an automorphism $\theta$ exists.
For each $r\in R$, $xr=\alpha(r)x$ and
$x^{-1}r=\alpha^{-1}(r)x^{-1}$. Applying $\theta$ to the first of
these, $x^{-1}\gamma(r)=\gamma\alpha(r)x^{-1}$, whence
$\alpha^{-1}\gamma=\gamma\alpha$ and
$\alpha^{-1}=\gamma\alpha\gamma^{-1}$. Thus $\alpha$ is
$\gamma$-reversible.

Conversely, suppose that $\alpha$ is $\gamma$-reversible and let
$\eta=\iota\gamma:R\rightarrow S$, where $\iota$ is the embedding of
$R$ in $S$.  The unit $x^{-1}$ in $S$ is such that
$x^{-1}\eta(r)=\eta\alpha(r)x^{-1}$. By the universal mapping
property for skew Laurent polynomial rings, as specified in
\cite[Exercise 1N]{GW}, there is a (unique) ring endomorphism
$\theta$ of $S$ such that $\theta|_R=\gamma$ and $\theta(x)=x^{-1}$.
Being self-inverse, $\theta\in \Aut(R)$.
\end{proof}

\begin{defn}
When $\alpha$ is $\gamma${\it-reversible} for some $\gamma\in
\Aut(R)$ such that $\gamma^2=\id_R$, we shall say that $S$ is a {\it
reversible} skew Laurent polynomial ring and that the automorphism
$\theta$ of $S$ such that $\theta(x)=x^{-1}$ and $\theta|_R=\gamma$
is the {\it reversing} automorphism of $S$ determined by $\gamma$.
\end{defn}

The main examples are the two pairs of related skew Laurent
polynomial rings discussed in the Introduction.
\begin{example}
\label{locenv} (i) Let $R=\F[y]$, let $\alpha,\gamma\in\Aut_\F(R)$
be such that $\alpha(y)=y+1$ and $\gamma(y)=-y$ and let
$S_1=R[x^{\pm 1};\alpha]$. Then $\alpha$ is $\gamma$-reversible,
$R[x;\alpha]=\F[y,x:xy-yx=x]$ is the enveloping algebra
$U(\mathfrak{g})$ of the two-dimensional non-abelian solvable Lie
algebra $\mathfrak{g}$ and $S_1$ is its localization, which we
denote $V(\mathfrak{g})$, at the powers of the normal element $x$.
The reversing automorphism $\theta_1$ of $V(\mathfrak{g})$
determined by $\gamma$ is such that $\theta_1(y)=-y$ and
$\theta_1(x)=x^{-1}$. It is well-known that $V(\mathfrak{g})$ is
simple if $\ch \F=0$, for example see \cite[Exercise 1V]{GW}.

(ii) Let $R=\F[y,t]$, let $\alpha,\gamma\in\Aut_\F(R)$ be such that
$\alpha(y)=y+1$, $\gamma(y)=-y$ and $\alpha(t)=t=\gamma(t)$. Let
$S=R[x^{\pm 1};\alpha]$. Then $\alpha$ is $\gamma$-reversible and
$R[x;\alpha]$ is the $\F$-algebra generated by $t,x$ and $y$ subject
to the relations
\[
xy-yx=tx,\quad tx=xt,\quad ty=yt.\] This is the {\it homogenized
enveloping algebra} $U_t(\mathfrak{g})$, with $\mathfrak{g}$ as in
(i), and $S$ is its localization, which we denote
$V_t(\mathfrak{g})$, at the powers of $x$. The reversing
automorphism $\theta$ of $V_t(\mathfrak{g})$ determined by $\gamma$
is such that $\theta(x)=x^{-1}$, $\theta(y)=-y$ and $\theta(t)=t$.
\end{example}

\begin{example} \label{qtorus}
(i) Let $q\in \F\backslash\{0\}$, let $R=\F[y^{\pm 1}]$ and let
$\alpha,\gamma\in\Aut_\F(R)$ be such that $\alpha(y)=qy$ and
$\gamma(y)=y^{-1}$. Then $\gamma^2=\id$ and $\alpha$ is
$\gamma$-reversible. Here the skew Laurent polynomial ring
$S=R[x^{\pm 1};\alpha]$ is the quantized coordinate ring
$W_q=\mathcal{O}_q((\F^*)^2)$ or, more informally, the quantum
torus, see \cite[p.16]{GW}. It is well-known that $W_q$ is simple if
$q$ is not a root of unity, for example see \cite[Corollary
1.18]{GW}. The reversing automorphism $\theta_q$ of $W_q$ determined
by $\gamma$ is such that $\theta_q(y)=y^{-1}$ and
$\theta_q(x)=x^{-1}$.

(ii) Let $R=\F[y^{\pm 1}, Q^{\pm 1}]$ and let $\alpha,\gamma\in
\Aut_\F(R)$ be such that $\alpha(y)=Qy$, $\alpha(Q)=Q=\gamma(Q)$ and
$\gamma(y)=y^{-1}$. The {\it generic quantum torus} is the skew
Laurent polynomial ring $W_Q=R[x^{\pm 1};\alpha]$. Then $xy=Qyx$,
$\alpha$ is $\gamma$-reversible and the reversing automorphism
$\theta$ of $W_Q$ determined by $\gamma$ is such that
$\theta(x)=x^{-1}$, $\theta(y)=y^{-1}$ and $\theta(Q)=Q$.
\end{example}

For the remainder of this section, let $R$ be a ring and let
$\alpha,\gamma\in \Aut(R)$ be such that $\gamma^2=\id_R$ and
$\alpha$ is $\gamma$-reversible. Let $S=R[x^{\pm 1};\alpha]$ and let
$\theta\in \Aut(S)$ be the reversing automorphism determined by
$\gamma$. We now identify some elements of $S^\theta$ and some
relations that hold between them.

\begin{lemma}
\label{srelations} For $r\in R$ and $n\geq 0$, let
$s_n(r):=rx^n+\gamma(r)x^{-n}$. In particular $s_0(r)=r+\gamma(r)$.
Then $s_n(r)\in S^\theta$. If $r,r^\prime\in R$ and
$rr^\prime=r^\prime r$ then
\begin{equation}
\label{rrprime}s_0(r)s_1(r^\prime)-s_1(r^\prime)s_0(\alpha^{-1}(r))=
s_1((\gamma(r)-\alpha^2(\gamma(r)))r^\prime).
\end{equation}
In particular, with $r^\prime=1$,
\begin{equation}
\label{sr}s_0(r)s_1(1)-s_1(1)s_0(\alpha^{-1}(r))=s_1(\gamma(r)-\alpha^2\gamma(r)).
\end{equation}
Also,
\begin{equation}
s_1(r)s_1(1)-s_1(1)s_1(\alpha^{-1}(r))=s_0(r-\alpha^{-2}(r)).
\label{sr2}
\end{equation}
\end{lemma}
\begin{proof} It is immediate from the definition of $\theta$ that $s_n(r)\in S^\theta$.
The relations \eqref{rrprime}, \eqref{sr} and \eqref{sr2} are
routinely checked using the equations
$\alpha^{-1}\gamma\alpha^{-1}=\gamma$ and
$\alpha\gamma\alpha^{-1}=\gamma\alpha^{-2}=\alpha^{2}\gamma$.
\end{proof}

\begin{prop}
The fixed ring $S^\theta$ is generated by the fixed ring $R^\gamma$
and the set $\{s_1(r):r\in R\}$.
 \label{gen1}
\end{prop}
\begin{proof}
Let $S_1$ be the subring of $S$ generated by $R^\gamma$ and
$\{s_1(r):r\in R\}$. It is clear that $S_1\subseteq S^\theta$. Let
$s=\sum_m^n r_ix^i\in S^\theta$, where each $r_i\in R$. Then
$s=\theta(s)=\sum_m^n \gamma(r_i)x^{-i}$ from which it follows that
$m=-n$, $r_0=\gamma(r_0)$ and, for $1\leq i\leq n$,
$r_{-i}=\gamma(r_i)$. Thus $s=r_0+\sum_1^n s_i(r_i).$ As $r_0\in
R^\gamma\subset S_1$, it now suffices to show that, for all $r\in R$
and all $i\geq 1$, $s_i(r)\in S_1$. This is certainly true when
$i=1$ and it follows inductively using the formula
\[
s_{i+1}(r)=s_i(r)s_1(1)-s_{i-1}(r).
\]
\end{proof}

The following result, whose hypothesis is satisfied if $R$ is left
Noetherian, by \cite[Corollary 26.13]{passman}, will be applicable
to give finite sets of generators for our examples.

\begin{prop}
Suppose that $R$ is finitely generated, as a left $R^\gamma$-module,
by $r_1=1,r_2,\ldots,r_n$. Then $S^\theta$ is generated by
$R^\gamma$ and $\{s_1(r_i):1\leq i\leq n\}$. \label{gen3}
\end{prop}
\begin{proof}
If $r=c_1r_1+c_2r_2+\ldots+c_nr_n$, where $c_1,c_2,\ldots,c_n\in
R^\gamma$, then $s_1(r)=c_1s_1(r_1)+c_2s_1(r_2)+\ldots+c_ns_1(r_n)$.
The result follows from Proposition~\ref{gen1}.
\end{proof}

\begin{cor}
\label{leacorhom} Let $S=V_t(\mathfrak{g})$ and $\theta$ be as in
Example~\ref{locenv}(ii) and suppose that $\ch \F\neq 2$. Then
$S^{\theta}$ is generated by $t$, $y^2$, $x+x^{-1}$ and
$yx-yx^{-1}$. \end{cor}
\begin{proof} Here $R^\gamma$ is generated by $t$ and $y^2$ and
$R=R^\gamma+R^\gamma y$ so the result follows from
Proposition~\ref{gen3}.
\end{proof}

\begin{cor}
\label{leacor} Let $S_1=V(\mathfrak{g})$ and $\theta_1$ be as in
Example~\ref{locenv}(i) and suppose that $\ch \F\neq 2$. Then
$S^{\theta_1}$ is generated by $y^2$ and $x+x^{-1}$.
\end{cor}
\begin{proof} Let $r=y^2/2$, $a_1=y^2=s_0(r)$,
$a_2=yx-yx^{-1}=s_1(y)$ and  $a_3=x+x^{-1}=s_1(1)$. Then $R^\gamma$
is generated by $a_1$ and $R=R^\gamma+R^\gamma y$ so, by
Proposition~\ref{gen3}, $ S_1^{\theta_1}$ is generated by $a_1, a_2$
and $a_3$. Also $\alpha^{-1}(r)=1$ and $\alpha^2\gamma(r)=r+2y+2$
so, by \eqref{sr}, $a_1a_3-a_3(a_1+1)=s_1(-2y-2)=-2a_2-2a_3,$
whence, as $\ch \F\neq 2$, $a_2$ is in the $\F$-subalgebra generated
by $a_1$ and $a_3$.
\end{proof}

\begin{cor}
\label{Qtcor} Let $S=W_Q$ and $\theta$ be as in
Example~\ref{qtorus}(ii). Then $S^\theta$ is generated by
$y+y^{-1}$, $Q$, $Q^{-1}$, $x+x^{-1}$ and $yx+y^{-1}x^{-1}$.
\end{cor}
\begin{proof} Here $R^\gamma$ is generated
by $Q$, $Q^{-1}$ and $y+y^{-1}$. For $n\geq0$, let $V_n$ be the
$\F[Q^{\pm 1}]$-submodule of $\F[y^{\pm 1},Q^{\pm 1}]$ generated by
$\{y^m:n\geq m\geq -n\}$. For $r,s\in \F[Q^{\pm1}]$,
$ry+sy^{-1}=s(y+y^{-1})+(r-s)y$, whence $V_1\subseteq
R^\gamma+R^\gamma y$. For $n\geq 2$,
\[ry^n+sy^{-n}=sy^n+sy^{-n}+((r-s)y^{n-1}+(r-s)y^{1-n})y+(s-r)y^{2-n},\]
so $V_n\subseteq R^\gamma+R^\gamma y+V_{n-1}.$ It follows,
inductively, that $V_n\subseteq R^\gamma+R^\gamma y$ for all $n$ and
hence that $R=R^\gamma+R^\gamma y$. The result follows from
Proposition~\ref{gen3}.
\end{proof}

\begin{cor}
\label{qtcor} Let $S_q=W_q$ and $\theta_q$ be as in
Example~\ref{qtorus}(i). Then $S_q^{\theta_q}$ is generated by
$y+y^{-1}$, $x+x^{-1}$ and $yx+y^{-1}x^{-1}$. If $q^2\neq 1$ then
$S^{\theta_q}$ is generated by $y+y^{-1}$ and $x+x^{-1}$.
\end{cor}
\begin{proof} Here $R^\gamma$ is generated
by $y+y^{-1}$ and, as in the proof of Corollary~\ref{Qtcor},
$R=R^\gamma+R^\gamma y$. The first conclusion follows from
Proposition~\ref{gen3}. Let $a_1=y+y^{-1}=s_0(y)$,
$a_2=x+x^{-1}=s_1(1)$ and $a_3=yx+y^{-1}x^{-1}=s_1(y)$. Note that
$s_0(\alpha^{-1}(y^{-1}))=qa_1$ and
$\gamma(y^{-1})-\alpha^2\gamma(y^{-1})=(1-q^2)y$. By \eqref{sr} with
$r=y^{-1}$, $ a_1a_2-qa_2a_1=(1-q^2)a_3$ so as $q^2\neq 1$, $a_3$ is
in the subalgebra generated by $a_1$ and $a_2$.
\end{proof}

\begin{rmk}
If $\ch \F=0$ in Example~\ref{locenv}(i) then $S_1^{\theta_1}$ is
simple as an easy consequence of \cite[Theorem 28.3(ii)]{passman}.
If $\ch \F\neq 2$ and $q$ is not a root of unity in
Example~\ref{qtorus}(i) then $S_q^{\theta_q}$ is simple for the same
reason. \label{simple}
\end{rmk}

\begin{rmk}
Let $\alpha$ be a $\gamma$-reversible automorphism of a ring $R$ and
let $S=R[x^{\pm 1};\alpha]$. Then $\alpha$ is also
$\alpha\gamma$-reversible and $\gamma\alpha$-reversible and there
are reversing automorphisms $\theta^\prime$ and
$\theta^{\prime\prime}$ of $S$ determined by $\alpha\gamma$ and
$\gamma\alpha$ respectively. Suppose that there is a
$\gamma$-reversible automorphism $\beta$ of $R$ such that
$\beta^2=\alpha^{-1}$. If $\ch \F\neq 2$ this is the case in
Example~\ref{locenv}(i), with $\beta(y)=y-\frac{1}{2}$, and if $q$
is a square in $\F$ it is the case in Example~\ref{qtorus}(ii), with
$\beta(y)=q^{-1/2}y$. Then $\beta$ extends to an automorphism of $S$
with $\beta(x)=x$  and
$\beta\theta^\prime\beta^{-1}=\theta=\beta^{-1}\theta^{\prime\prime}\beta$.
Thus the reversing automorphisms $\theta$, $\theta^\prime$ and
$\theta^{\prime\prime}$ are conjugate in $\Aut(S)$ and hence their
rings of invariants are isomorphic.
\end{rmk}

\section{Poisson algebras, quantization and deformation}
\label{PQD} In this section $\ch \F=0$. By a {\it Poisson algebra}
we mean a commutative $\F$-algebra $A$ with a bilinear product
$\{-,-\}:A\times A\rightarrow A$ such that $A$ is a Lie algebra
under $\{-,-\}$ and, for all $a\in A$, $\{a,-\}$ is an
$\F$-derivation of $A$. Such a product is a {\it Poisson bracket} on
$A$. For $a\in A$, the derivation $\{a,-\}$ is a {\it hamiltonian}
derivation (or hamiltonian vector field) and is written $\ham a$.

A subalgebra $B$ of $A$ is a {\it Poisson subalgebra} of $A$ if
$\{b,c\}\in B$ for all $b,c\in B$ and an ideal $I$ of $A$ is a {\it
Poisson ideal} if $\{i,a\}\in I$ for all $i\in I$ and all $a\in A$.
If $I$ is a Poisson ideal of $A$ then $A/I$ is a Poisson algebra in
the obvious way: $\{a+I,b+I\}=\{a,b\}+I$.  The {\it Poisson centre}
of a Poisson algebra $A$ is $\PZ(A):=\{a\in A: \{a,b\}=0\mbox{ for
all }b\in A\}$.

An $\F$-automorphism $\pi$ of a Poisson algebra $A$ is a {\it
Poisson automorphism} of $A$ if $\{\pi(a),\pi(b)\}=\pi(\{a,b\})$ for
all $a,b\in A$, in which case $\pi^{-1}$ is also a Poisson
automorphism. If $\pi$ is a Poisson automorphism of $A$ then the
ring of invariants $A^\pi$ is a Poisson subalgebra of $A$.

\begin{defns}
\label{Poissonind} Let $T$ be an $\F$-algebra with a central
non-unit non-zero-divisor $t$ such that $B:=T/tT$ is commutative.
Then $[-,-]$ in $T$ induces a well-defined Poisson bracket $\{-,-\}$
on $B$ by the rule
\[\{\ov u,\ov v\}=\ov{t^{-1}[u,v]}\mbox{ for all }
\ov u=u+tT,\ov v=v+tT\in B.\] For more detail, see \cite[Chapter
III.5]{BGl}. Following \cite{BGl}, we say that $T$ is a {\it
quantization} of the Poisson algebra $B$. Let $\lambda\in \F$ be
such that the central element $t-\lambda$ is a non-unit in $T$ and
let $T_\lambda=T/(t-\lambda)T$. We shall refer to $T_\lambda$ as a
{\it deformation} of $B$.

Now suppose that there exists $\theta\in \Aut_\F (T)$ such that
$\theta(t)=t$. Then there is a well-defined $\F$-automorphism $\pi$
of $B$ such that $\pi(\ov u)=\ov{\theta(u)}$ for all $\ov u\in B$.
Let $a=u+tT, b=v+tT\in B$. Then
\begin{align*}
\{\pi(a),\pi(b)\}&=\{\theta(u)+tT,\theta(v)+tT\}=
t^{-1}[\theta(u),\theta(v)]+tT\\
&=\theta(t^{-1}[u,v])+tT=\pi(\{u+tT,v+tT\})=\pi(\{a,b\}).
\end{align*}
Thus $\pi$ is a Poisson automorphism of $B$. There is also an
automorphism $\theta_\lambda$ of the deformation $T/(t-\lambda)T$
with $\theta_\lambda(\ov u)=\ov{\theta(u)}$ for all $\ov u\in
T/(t-\lambda)T$. We shall refer to $\theta$ as a quantization of
$\pi$ and to $\theta_\lambda$ as a {\it deformation} of $\pi$.

In this situation, $t\in T^{\theta}$, where it must be a regular
central non-unit, so the $\F$-algebra $C:=T^{\theta}/tT^{\theta}$
becomes a Poisson algebra. On the other hand the ring of invariants
$B^{\pi}$ is a Poisson subalgebra of $B$. There is an injective
$\F$-algebra homomorphism $\psi:C\rightarrow B^{{\pi}}$ given by
$\psi(u+tT^{\theta})=u+tT$. If $\theta$ has the property that
$\{u\in T: u-\theta(u)\in tT\}= tT+T^{\theta}$ then $\psi$ is an
isomorphism and we may identify the Poisson algebras $C$ and
$B^{\pi}$.
\end{defns}

\begin{example}\label{Plocenv}
In Example~\ref{locenv}, $V_t(\mathfrak{g})$ is a quantization of
the Poisson algebra $B:=\F[x^{\pm 1},y]$, with the Poisson bracket
such that $\{x,y\}=x$, and $V(\mathfrak{g})$ is a deformation of
$B$. By the quotient rule for the derivation $\ham y$, $
\{x^{-1},y\}=-x^{-1}$. The reversing automorphisms $\theta$ and
$\theta_1$ of $V_t(\mathfrak{g})$ and $V(\mathfrak{g})$ are,
respectively, a quantization and a deformation of the Poisson
automorphism $\pi$ of $B$ such that $\pi(x)=x^{-1}$ and $\pi(y)=-y$.
There are three obvious invariants under $\pi$, namely $a_1:=y^2$,
$a_2:=y(x-x^{-1})$  and $a_3:=x+x^{-1}$. It is a routine matter to
check that these generate $B^\pi$ and that $a_2^2=a_1(a_3^2-4)$.
Thus, if $A=\F[x_1,x_2,x_3]$ and $f=x_1(4-x_3^2)+x_2^2$, then there
is a surjection $\phi: A\rightarrow B^\pi$ with $f\in \ker \phi$. As
$B^\pi$ and $A/fA$ are domains of Krull dimension $2$, they are
isomorphic.
\end{example}

\begin{example} \label{Pqtorus}
In Example~\ref{qtorus}, let $t=Q-1$, which is a regular central
non-unit in $W_Q$, and let
 $B=W_Q/tW_Q$ which we identify with $\F[x^{\pm 1},y^{\pm 1}]$.
 Then $W_Q$ is a quantization of the Poisson algebra $B$, where the Poisson bracket
 is such that $\{x,y\}=xy$, and, taking $\lambda=q-1$, $W_q$ is a deformation.
 In $B$, $\ham x=xy\partial/\partial y$ and $\ham y=-xy\partial/\partial x$.
Consequently $\{x,y^{\pm 1}\}=-xy^{\pm 1}$, $\{x^{\pm 1},y\}=-x^{\pm
1}y$ and $\{x^{\pm 1},y^{\pm 1}\}=x^{\pm 1}y^{\pm 1}$. There is a
Poisson automorphism $\pi$ of $B$ such that $\pi(x)=x^{-1}$ and
$\pi(y)=y^{-1}$ and the reversing automorphisms $\theta$, of $W_Q$,
and $\theta_q$, of $W_q$ are, respectively, a quantization and
deformation of $\pi$.

The ring of invariants of the commutative Laurent polynomial
$\Z[y^{\pm 1},x^{\pm 1}]$ for the automorphism $x\mapsto x^{-1}$,
$y\mapsto y^{-1}$ is discussed in \cite[Example 3.5]{lor}. The ring
of invariants is generated by $a_1:=y+y^{-1}, a_2:=x+x^{-1}$ and
$a_3:=xy+x^{-1}y^{-1}$ and is isomorphic to
$\Z[x_1,x_2,x_3]/f\Z[x_1,x_2,x_3]$ where
$f=x_1x_2x_3-x_1^2-x_2^2-x_3^2+4$. The same calculations are valid
over $\F$, where $B^\pi \simeq \F[x_1,x_2,x_3]/f\F[x_1,x_2,x_3]$.
\end{example}

\begin{rmk}
\label{invexactP} Let $A=\F[x_1,x_2,x_3]$. In each of the
Examples~\ref{Plocenv} and \ref{Pqtorus}, we have identified
$B^\pi$ as the coordinate ring $A/fA$ of a surface determined by an
irreducible polynomial $f$. In each case there is a Poisson bracket
on $A$, determined by $f$ as in the following definition, such that
$B^\pi$ is a factor of $A$ as a Poisson algebra.
\end{rmk}

\begin{defn}
Let $A=\F[x_1,x_2,x_3]$ and let $f\in A$. There is a Poisson bracket
$\{-,-\}_f$ on $A$ given by $\{x_1,x_2\}_f=\partial f/\partial x_3$,
$\{x_2,x_3\}_f=\partial f/\partial x_1$ and $\{x_3,x_1\}_f=\partial
f/\partial x_2$. Such brackets are considered, for example, in
\cite{pich}, \cite[p.1312 (1) with $n=3$ and $\lambda=1$]{Odesskii+}
and \cite[p. 252]{dml}.

For $g,h\in A$, \[\{g,h\}_f=\left|
\begin{array}{ccc}
  \frac{\partial f}{\partial x_1}& \frac{\partial f}{\partial x_2} &
  \frac{\partial f}{\partial x_3} \\
  \frac{\partial g}{\partial x_1}& \frac{\partial g}{\partial x_2} &
  \frac{\partial g}{\partial x_3} \\
  \frac{\partial h}{\partial x_1}& \frac{\partial h}{\partial x_2} &
  \frac{\partial h}{\partial x_3}\\
\end{array}
\right|\] from which it is clear that $f\in \PZ(A)$ and hence that
$fA$ is a Poisson ideal of $A$. Writing $f_i=\partial f/\partial
x_i$ for $i=1,2,3$, the hamiltonian derivations of $A$ for the three
generators are
\begin{align*}
\ham x_1&=f_3\partial/\partial x_2-f_2\partial/\partial x_3;\\ \ham
x_2&=f_1\partial/\partial x_3-f_3\partial/\partial x_1;\\ \ham
x_3&=f_2\partial/\partial x_1-f_1\partial/\partial x_2.
\end{align*}

We shall call a Poisson bracket on $A$ {\it exact} (determined by
$f$) if it has the form $\{-,-\}_f$ for some $f\in A$.

\label{Poissonf}
\end{defn}

\begin{examples} \label{poissoninv} In \ref{Plocenv},
it is a straightforward exercise to check that
\[\{a_1,a_2\}=-2a_1a_3,\,
\{a_2,a_3\}=4-a_3^2\mbox{ and }\{a_3,a_1\}=2a_2.
\]
Thus, as a Poisson algebra, $B^\pi$ is the factor $A/fA$ of $A$
under the exact Poisson bracket $\{-,-\}_f$, where
$f=x_1(4-x_3^2)+x_2^2$ and each $a_i=x_i+fA$.

Similarly, in Example~\ref{Pqtorus}, it can be checked that
\[\{a_1,a_2\}=a_1a_2-2a_3,\, \{a_2,a_3\}=a_2a_3-2a_1\mbox{ and }
\{a_3,a_1\}=a_1a_3-2a_2
\]
so that the Poisson algebra $B^\pi$ is again a Poisson factor
$A/fA$ of $A$, with the exact Poisson bracket $\{-,-\}_f$, where, on
this occasion, $f=x_1x_2x_3-x_1^2-x_2^2-x_3^2+4$.
\end{examples}

\begin{rmk}
In Example~\ref{locenv}(ii) the automorphism $\alpha$, used in
constructing $S=V_t(\mathfrak{g})$ as a skew Laurent polynomial
ring, extends to a $\F$-automorphism of $S$ such that $\alpha(x)=x$,
$\alpha(t)=t$ and $\alpha(y)=y+t\equiv y \bmod tS$. Thus $\alpha$
quantizes the identity automorphism on $B=\F[y,x^{\pm 1}]$. The same
is true in Example~\ref{qtorus}(ii) where the automorphism $\alpha$
of $W_Q$ such that $\alpha(x)=x, \alpha(y)=Qy$ and $\alpha(Q)=Q$
quantizes the identity automorphism on $B=\F[y^{\pm 1},x^{\pm 1}]$.
\end{rmk}

\begin{rmk}
\label{simPoiss} In Examples~\ref{Plocenv} and \ref{Pqtorus}, the
Poisson algebra $B$ is simple. In each case $B$ is a localization of
$\F[x,y]$ and the hamiltonian derivations $\{x,-\}$ and $\{y,-\}$
have the forms $u\partial/\partial y$ and $v\partial/\partial x$
respectively, where $u$ and $v$ are units in $B$. Consequently any
non-zero Poisson ideal of $B$ must intersect $\F[x,y]$ in a non-zero
ideal invariant under the derivations $\partial/\partial x$ and
$\partial/\partial y$ and therefore cannot be proper.

There is an example in \cite{alevfarkas} of a Poisson automorphism
($x\mapsto -x, y\mapsto -y$) of a simple Poisson algebra ($\F[x,y]$
with $\{x,y\}=1$) such that the ring of invariants is not a simple
Poisson algebra, having a maximal ideal (generated by $x^2,xy,y^2$)
that is Poisson. The situation for our two examples is similar in
that, although $B$ is simple, there are finitely many maximal ideals
 of $B^{\pi}$ that are Poisson. In
Example~\ref{Plocenv}, these are the two maximal ideals
$(a_1,a_2,a_3\pm 2)$  while, in Example~\ref{Pqtorus}, they are
 $(a_1-2,a_2-2,a_3-2)$, $(a_1+2,a_2+2,a_3-2)$,
$(a_1+2,a_2-2,a_3+2)$ and $(a_1-2,a_2+2,a_3+2)$. Thus $B^\pi$ is
not simple in either example.
\end{rmk}

\section{Filtrations}
\label{filtrations} We have identified finite sets of generators for
the rings of invariants for the four reversing automorphisms in
Examples~\ref{locenv} and \ref{qtorus}. Identifying finite sets of
relations is much more technical. A model is given by the argument,
using Krull dimension and domain recognition, sketched in
Examples~\ref{poissoninv} for the invariants of the automorphism
$\pi$ of $\F[x^{\pm 1},y]$ with $\pi(x)=x^{-1}$ and $\pi(y)=-y$.
Essentially, we need to impose sufficiently many relations, from
those identified in Lemma~\ref{srelations}, to obtain a domain of
the correct dimension. The appropriate dimension is Gelfand-Kirillov
dimension, for which we refer to \cite{KL} and \cite{McCR}. Other
key methods involve Bergman's Diamond Lemma, for which we refer to
\cite[Appendix I.11]{BGl}, and filtrations with their associated
graded rings, for which references are \cite[\S 1.6]{McCR} and
\cite[Appendix I.12]{BGl}. The methods will occasionally be
sensitive to the choice of filtration. The filtrations considered
are slightly more general than those described in
\cite[I.12.2(c)]{BGl}.

The orderings on monomials that we use in applying the Diamond Lemma
are
modifications of the length-lexicographic ordering. Let $n\geq 1$,
let $M_n$ be the free monoid on $\{ z_1, z_2, \ldots, z_n\}$ and let
$F_n$ be the free algebra $\F\langle z_1, z_2, \ldots,
z_n\rangle$.  By a {\it degree function} on $M_n$, we mean a monoid
homomorphism $d:M_n\rightarrow (\mathbb{N}_0,+)$. Such a function is
determined by its values on $z_1, z_2, \ldots, z_n$.
\begin{defns}
\label{dfiltrations} Given a degree function $d:M_n\rightarrow
(\mathbb{N}_0,+)$ we modify the length-lexicographic ordering
$\preceq_{\lex}$ by ordering words first by the degree function $d$
and then lexicographically with $z_1\succ z_2\succ\ldots\succ z_n$.
More formally, we define $m\leq m^\prime$ if and only if either
$d(m)<d(m^\prime)$ or $d(m)=d(m^\prime)$ and $m\preceq_{\lex}
m^\prime$. We do not require that $d(z_i)\geq d(z_j)$ whenever
$i\leq j$ so we may have $z_2>z_1$ although $z_1\succ z_2$. We shall
refer to this as the {\it{ $d$-length-lexicographic ordering }}. It
is clearly a semigroup ordering and if $d(z_i)>0$ for each $i$ it
has DCC.

If $d(z_i)=0$ for some $i$ then a semigroup ordering with DCC can be
defined using a complementary degree function $e$, such that
$e(z_i)=0$ if $d(z_i)>0$ and $e(z_i)=1$ if $d(z_i)=0$. The ordering
is then given by the rules: $m\leq m^\prime$ if and only if either
$d(m)<d(m^\prime)$ or ($d(m)=d(m^\prime)$ and $e(m)<e(m^\prime)$) or
($d(m)=d(m^\prime)$ and $e(m)=e(m^\prime)$ and $m\preceq_{\lex}
m^\prime$). We shall refer to this as the {\it augmented
$d$-length-lexicographic ordering.} It has DCC because for each
$m\in M_n$ there are only finitely many monomials $<m$.
\end{defns}

Let $A$ be an algebra with a presentation of the form $F_n/I$, where
$I$ is an ideal of $F_n$, and let $x_i=z_i+I\in A$, $1\leq i\leq n$.
Let $d:M_n\rightarrow (\mathbb{N}_0,+)$ be a degree function. Set
$A_0=\F$, and, for $i\geq 1$, let $A_i$ be the $\F$-subspace of $A$
spanned by the images $m+I$ in $A$ of words $m\in M_n$ with
$d(m)\leq i$. Then $A_0\subseteq A_1\subseteq A_2\ldots $ is a
filtration of $A$. We shall call this the {\it $d$-standard
filtration} of $A$.

Suppose now that $I$ is the ideal generated by the elements
$w_\sigma-f_\sigma$, for some reduction system
$S=\{(w_\sigma,f_\sigma)\}$. Let $\leq$ be a semigroup ordering on
$M_n$ that has DCC and is compatible with $S$. We say that a degree
function $d:M_n\rightarrow \mathbb{N}_0$ is {\it compatible} with $S$
if, for each $(w_\sigma,f_\sigma)\in S$, $f_\sigma$ is a linear
combination of words $m$ with $d(m)\leq d(w_\sigma)$. If $\leq$ is
the augmented $d$-length-lexicographic ordering then compatibility
of $d$ with $S$ is a consequence of compatibility of $\leq$ with
$S$.

The following Proposition will be applicable to identify associated
graded rings for the filtrations that we use.
\begin{prop}
\label{grbasis} Let $A$, $M_n$, $F_n$ and $S$ be as above. Let
$\leq$ be a semigroup ordering on $M_n$, with DCC, that is
compatible with $S$ and let $d$ be a degree function that is
compatible with $S$. Denote by $\Irr(M_n)$ the set of images $m+I$
in $A$ of those words $m$ in $M_n$ that are irreducible with respect
to $S$. Suppose that all ambiguities in $S$ are resolvable. Then,
for the $d$-standard filtration, $\gr(A)$ has basis $\{\ov{m+I}:m\in
\Irr(M_n)\}$ and there is a vector space isomorphism
$\psi:A\rightarrow \gr (A)$ given by $\psi(m+I)=\ov{m+I}$ for all
$m+I\in \Irr(M_n)$.
\end{prop}
\begin{proof} By the Diamond Lemma, $\Irr(M_n)$ is
a basis for $A$. Set $B_0=\F$, and, for $i\geq 1$, let $B_i$ be the
$\F$-subspace spanned by the images $m+I$ in $A$ of the irreducible
words $m$ with $d(m)\leq i$. Thus $B_i\subseteq A_i$ and $B_i$ has
basis $B_i\cap \Irr(M_n)$. We claim that $B_i=A_i$ for each $i$.
Suppose not. By DCC, there exists a word $m\in M_n$ that is minimal,
under $\leq$, with the property that $m+I\in A_{d(m)}\backslash
B_{d(m)}$. Then $m$ cannot be irreducible so there exists
$(w_\sigma, f_\sigma)\in S$ such that $m=aw_\sigma b$ for some
$a,b\in M_n$. Then $m=af_\sigma b$ is a linear combination of words
$w$ with $w<m$ and $d(w)\leq d(m)$, whence $m+I\in B_{d(m)}$,
contradicting the minimality of $m$. Therefore, for all $i$,
$A_i=B_i$, $A_i$ has basis $A_i\cap \Irr(M_n)$ and, in $\gr(A)$,
each summand $A_i/A_{i-1}$ has a basis consisting of the elements
$\ov{m+I}$ where $m+I\in (A_i\cap \Irr(M_n))\backslash A_{i-1}$.
Therefore $\{\ov{m+I}:m\in \Irr(M_n)\}$ is a basis of $\gr(A)$ and
there is a vector space isomorphism $\psi:A\rightarrow \gr(A)$ given
by $\psi(m+I)=\ov{m+I}$ for all $m+I\in \Irr(M_n)$.
\end{proof}

\section{Invariants for the localized enveloping algebra and its homogenization}
\label{invlocenv} Having identified, in Section \ref{revbasics},
generators for the rings of invariants of our principal examples of
reversing automorphisms, we now aim to identify full sets of
defining relations, beginning in this section with those specified
in Example~\ref{locenv}. We shall assume that $\ch \F\neq 2$. The
methods for $V(\mathfrak{g})$ and $V_t(\mathfrak{g})$ are similar
and we begin with the latter.  We know from
Corollary~\ref{leacorhom} that $V_t(\mathfrak{g})^{\theta}$ is
generated by $t$, $y^2$, $x+x^{-1}$ and $yx-yx^{-1}$.

\begin{prop}
\label{invhomlocenv}Let $S=V_t(\mathfrak{g})$ and the
reversing automorphism $\theta$ be as in Example~\ref{locenv}(ii). Let
$T$ be the $\F$-algebra generated by $t, x_1, x_3$ and
$x_2$ subject
to the relations \begin{align}\label{tcent} x_it&=tx_i\mbox{ for }i=1,2,3,\\
\label{X12rel} x_1x_2&=x_2x_1-2tx_3x_1+3t^2x_2+2t^3x_3,\\
x_2x_3&=x_3x_2-tx_3^2+4t,\label{X23rel}\\
x_1x_3&=x_3x_1-t^2x_3-2tx_2. \label{X13rel}
\end{align}

(i) $T$ is an iterated skew polynomial ring
$\F[t,x_3][x_2;\delta][x_1;\sigma,\delta_1]$, where $\delta$ is a
derivation of $\F[t,x_3]$, $\sigma$ is an automorphism of
$\F[t,x_3][x_2;\delta]$ and $\delta_1$ is a $\sigma$-derivation of
$\F[t,x_3][x_2;\delta]$.

(ii) Let $g=(4-x_3^2)x_1+x_2^2+3tx_3x_2+t^2x_3^2+4t^2$. Then $g$ is
a central element of $T$ and $T/gT$ is
a domain.

(iii) $S^{\theta}$ is isomorphic to $T/gT$.
\end{prop}
\begin{proof}
(i) Let $R_1$ be the commutative polynomial ring $\F[t,x_3]$ and let
$\delta$ be
 the $\F$-derivation $t(4-x_3^2)\partial/\partial x_3$ of $R_1$.
 Let $R_2=R_1[x_2;\delta]$
 so that \eqref{X23rel} is satisfied.
 Let $F$ be the free algebra
 $\F\langle u, z_2, z_3\rangle$ and $I$ be the ideal of $F$ generated by
 $g:=uz_2-z_2u, h:=uz_3-z_3u$
 and $f:=
 z_2z_3-z_3
 z_2+uz_3^2-4u$. Then, by \cite[Proposition 1]{fddaj}, there is an isomorphism
 $\phi:R_2\rightarrow F/I$ with $\phi(t)=\ov u,
 \phi(x_2)=\ov z_2$ and
 $\phi(x_3)=\ov z_3$.

 There exists $\tau\in \Aut_\F(F)$ such
 that $\tau(u)=u$, $\tau(z_3)=z_3$ and $\tau(z_2)=z_2-2uz_3$,
 its inverse being such
 that $\tau^{-1}(u)=u$, $\tau^{-1}(z_3)=z_3$ and
 $\tau^{-1}(z_2)=z_2+2uz_3$. It is readily
 checked that $\tau(g)=g-2uh$, $\tau(h)=h$ and $\tau(f)=f-2hz_3$, whence $\tau(I)=I$ and
 there is an induced $\F$-automorphism $\sigma$ of $R_2$ such that
$\sigma(t)=t$, $\sigma(x_3)=x_3$ and $\sigma(x_2)=x_2-2tx_3$.

A left $\tau$-derivation $\partial$ of $F$ is determined by
specifying $\partial(u)$, $\partial(z_3)$ and $\partial (z_2)$ and
using the definition of a left $\tau$-derivation \cite[p. 33]{GW} to
extend to arbitrary elements of $F$. Here we set $\partial(u)=0$,
$\partial(z_3)=-u^2z_3-2uz_2$ and $\partial(z_2)=2u^3z_3+3u^2z_2$.
Modulo $I$, $u$ commutes with $z_3$ and $z_2$ so we find that
\begin{align*}
\partial(z_2z_3)&\equiv 2u^2z_2z_3+4u^2z_3z_2 +4u^3z_3^2-2uz_2^2\bmod I,\\
\partial(z_3z_2)&\equiv 2u^2z_3z_2 +2u^3z_3^2-2uz_2^2\bmod I\mbox{ and}\\
\partial(uz_3^2)&\equiv-2u^2z_2z_3-2u^2z_3z_2-2u^3z_3^2\bmod I,\end{align*}
 whence $\partial(f)\in I$.
Also $\partial(g)\in I$ and $\partial(h)\in I$. Thus
$\partial(I)\subseteq I$ so
 there is an induced $\sigma$-derivation $\delta_1$ of $R_2$ such that
$\delta_1(t)=0$, $\delta_1(x_3)=-t^2x_3-2tx_2$ and
$\delta_1(x_2)=2t^3x_3+3t^2x_2$. By \cite[Proposition 1]{fddaj},
$R_2[x_1;\sigma,\delta_1]$ is the $\F$-algebra generated by $t, x_1,
x_2$ and $x_3$ subject to the relations \eqref{tcent},
\eqref{X12rel}, \eqref{X23rel} and \eqref{X13rel}.

(ii) Let $\ad_1$ (\emph{resp} $\ad_3$) be the inner derivation of
$T$ given by $a\mapsto x_1a-ax_1$ (\emph{resp} $a\mapsto
x_3a-ax_3$). Then, using \eqref{X13rel}, \eqref{X12rel} and
\eqref{X23rel},
\begin{align*}
\ad_1(g)&=\ad_1(x_2^2)+\ad_1(x_3^2)(t^2-x_1)+3t\ad_1(x_3x_2)\\
&=-4tx_3x_2x_1+6t^2x_3^2x_1-8t^2x_1+6t^2x_2^2-2t^3x_3x_2-6t^4x_3^2+8t^4\\&\quad
-4t^3x_3x_2-8t^4
+4tx_3x_2x_1+8t^2x_1-6t^2x_3^2x_1+6t^3x_3x_2+6t^4x_3^2-6t^2x_2^2\\
&=0
\end{align*}
and
\begin{align*}
\ad_3(g)
&=\ad_3(x_2^2)+(4-x_3^2)\ad_3(x_1)+3x_3\ad_3(x_2)\\
&=2tx_3^2x_2-2t^2x_3^3-8tx_2+8t^2x_3+4t^2x_3+8tx_2\\
&\quad -t^2x_3^3-2tx_3^2x_2+3t^2x_3^3-12t^2x_3\\
&=0.
\end{align*}
Thus $gx_1=x_1g$ and $gx_3=x_3g$. By \eqref{X13rel} and
\eqref{tcent}, $g$ commutes with $tx_2$. As $t$ is central and $T$
is a domain, $g$ commutes with $x_2$. Therefore $g$ is central.

Let
 $d=4-x_3^2$ and $e=x_2^2+t^2x_3^2+3tx_3x_2+4t^2$.
 Applying \cite[Proposition 1]{normal} to the central element $g=dx_1+e$
in the ring $T=R_2[x_1;\sigma,\delta_1]$, we see that $dR_2$ is an
ideal of $R_2$ and that, provided $e$ is regular modulo $dR_2$,
$T/gT$ is a domain. Now $R_2/dR_2$ is commutative and may be
identified with $C/(4-x_3^2)C$ where $C$ is the commutative
polynomial ring $\F[t,x_2,x_3]$. Thus $R_2/dR_2$ has two minimal
primes $P_1$ and $P_2$, generated by the images of $2-x_3$ and
$2+x_3$, and intersecting in $0$. The set of zero-divisors in
$R_2/dR_2$ is $P_1\cup P_2$ and $e+dR_2\notin P_1\cup P_2$. Hence
$T/gT$ is a domain and $gT$ is a (completely) prime ideal of $T$.

(iii) We have seen in Corollary~\ref{leacorhom} that, in the
notation of Lemma~\ref{srelations}, $S^{\theta}$ is generated by
$t$, $a_1=y^2=s_0(\frac{1}{2}y^2)$, $a_2=yx-yx^{-1}=s_1(y)$ and
$a_3=x+x^{-1}=s_1(1)$. We now check that \eqref{tcent},
\eqref{X12rel}, \eqref{X23rel} and \eqref{X13rel} hold when $x_1,
x_2$ and $x_3$ are replaced by $a_1, a_2$ and $a_3$.
 This is certainly true for
\eqref{tcent}, $t$ being central in $S$.

Note that $s_0(y^i)=0$ if $i$ is odd, that $s_0(1)=2$ and that
$s_0(\alpha^{-1}(\frac{1}{2}y^2))=a_1+t^2$. By \eqref{sr},
\[a_1a_3-a_3(a_1+t^2)=s_1(-2ty-2t^2)=-2ta_2-2t^2a_3,\] whence
$a_1, a_2$ and $a_3$ satisfy \eqref{X13rel}, and, by \eqref{sr2},
\[a_2a_3=s_1(1)s_1(y-t)+s_0(2t)=a_3a_2-ta_3^2+4t,\] whence
$a_1, a_2$ and $a_3$ satisfy \eqref{X23rel}.
 For \eqref{X12rel}, note that, by \eqref{rrprime},
 \begin{equation}
 \label{a1a2}a_1a_2=a_2(a_1+t^2)+s_1(-2ty^2-2t^2y)=a_2a_1+t^2a_2+s_1(-2ty^2)-2t^2a_2.
 \end{equation}
 Let $r=-\frac{1}{3}y^3-ty^2-\frac{2}{3}t^2y$. Then
 $s_0(r)=s_0(-ty^2)=-2ta_1$,
 $s_0(\alpha^{-1}(r))=s_0(-\frac{1}{3}y^3+\frac{1}{3}t^2y)=0$ and
 $\gamma(r)-\alpha^2\gamma(r)=-2t^2y$.
By
 \eqref{sr} and \eqref{X13rel} for $a_1,a_2,a_3$,
\[s_1(-2ty^2)=-2ta_1a_3=-2ta_3a_1+2t^3a_3+4t^2a_2.\] Combining this
with \eqref{a1a2} shows that $a_1, a_2$ and $a_3$ satisfy
\eqref{X12rel}.

Thus \eqref{tcent}, \eqref{X12rel}, \eqref{X23rel} and
\eqref{X13rel} hold when $x_1, x_2$ and $x_3$ are replaced by $a_1,
a_2$ and $a_3$ so there is a surjective ring homomorphism
$\eta:T\rightarrow S^{\theta}$ such that $\eta(t)=t, \eta(x_1)=a_1,
\eta(x_2)=a_2$ and $\eta(x_3)=a_3$. It remains to show that $\ker
\eta=gT$.

A simple calculation, using \eqref{X13rel}, \eqref{X12rel} and
\eqref{X23rel}, shows that
\begin{equation}
\label{uv2} a_1a_3^2=a_3^2a_1-4ta_3a_2-8t^2.
\end{equation}
In $S^{\theta}$,
\begin{align*}
a_2^2&=(yx-yx^{-1})^2\\
&=y^2(x^2-2+x^{-2})+yt(x^2-x^{-2})\\
&=a_1(a_3^2-4)+ta_2a_3\\
&=(a_3^2-4)a_1-4ta_3a_2-8t^2+ta_3a_2-t^2a_3^2+4t^2\mbox{ (by \eqref{uv2} and \eqref{X23rel})}\\
&=(a_3^2-4)a_1-3ta_3a_2-4t^2-t^2a_3^2.
\end{align*}
Therefore $g\in \ker \eta$.

To show that $gT=\ker \eta$, we shall use Gelfand-Kirillov
dimension and a filtration of $T$. Let $F$ be the free
algebra $F_4$ and $M$ be the free
monoid  $M_4$. It will be convenient to write
$z_4$ as $u$. Let $\psi:F\rightarrow T$ be the surjective
homomorphism such that $\psi(z_i)=x_i$, $1\leq i\leq 3$, and
$\psi(u)=t$. Thus we may identify $T$ and $F/\ker \psi$.
Let $d$ be the degree function such that
\begin{equation}
d(z_1)=6,\;d(z_2)=4 ,\;  d(z_3)=2\mbox{ and } d(u)=1.
\end{equation}
Consider the $d$-standard filtration
 on $T$ and note that, as $T$ has been presented as an iterated skew
polynomial ring $\F[t,x_3][x_2;\delta][x_1;\sigma,\delta_1]$, with
$\sigma$ an automorphism, it has a basis
$\{t^ix_3^jx_2^kx_1^\ell\}$. It follows from this that, in
presentation of $T$ given in (i), and with the $d$-length
lexicographic ordering, all ambiguities are resolvable.
Alternatively, this may be checked directly. Computing degrees of
the monomials appearing in \eqref{X12rel}, \eqref{X23rel}, and
\eqref{X13rel} and applying Proposition~\ref{grbasis}, we see that
$\gr(T)$ is a commutative polynomial ring in four variables so, by
\cite[Proposition 6.6]{KL}, $\GK (T)=4$. On the other hand, it
follows from \cite[Proposition 3.5]{KL} that
$\GK(S)=\GK(\F[x^{\pm1},t][y;-tx\partial/\partial x])=3$. By
\cite[Corollary 26.13(ii)]{passman}, $S$ is finitely generated as a
right module over $S^{\theta}$, so, by \cite[Proposition
8.2.9]{McCR}, $\GK(S^{\theta})=3$.

As $S$ is a domain, so too is its subalgebra $S^{\theta}$. Therefore
$\ker \eta$ is a prime ideal $P$, say, of $T$, such that $T/P\simeq
S^{\theta}$. By \cite[Corollary 3.16]{KL},
\[4=\GK(T)\geq\GK(T/P)+\hgt(P)=\GK(S^{\theta})+\hgt(P)=3+\hgt(P).\]
Hence $\hgt(P)\leq 1$. As $T$ and $T/gT$ are domains, by (i) and
(ii), and as $0\neq g\in P$, it must be the case that
$gT=P=\ker\eta$. Therefore $S^{\theta}\simeq T/gT$.
\end{proof}

\begin{rmk}
Having identified invariants for the quantization $\theta$ of $\pi$,
we proceed to consider the deformations. All the deformations
$S/(t-\lambda)S$, $\lambda\in \F\backslash\{0\}$, are isomorphic and
the same is true of $T/(t-\lambda)T$. So we shall only consider the
case $\lambda=1$.
\end{rmk}

\begin{prop}
\label{thinvlocenv}Let $S_1=V(\mathfrak{g})$ and let the reversing
automorphism $\theta_1$ be as in Example~\ref{locenv}(i). Suppose
that $\ch \F\neq 2$. Let $T_1$ be the $\F$-algebra generated by
$x_1, x_2$ and $x_3$ subject to the relations
\begin{align}
x_1x_3&=x_3x_1-x_3-2x_2, \label{X12not}\\
x_2x_3&=x_3x_2-x_3^2+4,\label{X23not}\\
\label{X13not} x_1x_2&=x_2x_1-2x_3x_1+3x_2+2x_3.
\end{align}
(i) $T_1$ is an iterated skew polynomial ring
$\F[x_3][x_2;\delta][x_1;\sigma,\delta_1]$, where $\delta$ is a
derivation of $\F[x_3]$, $\sigma$ is an automorphism of
$\F[x_3][x_2;\delta]$ and $\delta_1$ is a $\sigma$-derivation of
$\F[x_3][x_2;\delta]$.

(ii) Let $p=(4-x_3^2)x_1+x_2^2+3x_3x_2+x_3^2+4$. Then $p$ is a
central element of $T_1$ and $T_1/pT_1$ is a domain.

(iii) $S_1^{\theta_1}$ is isomorphic to $T_1/pT_1$.
\end{prop}
\begin{proof}
This can be proved by the same methods as
Proposition~\ref{invhomlocenv}. The details are somewhat simpler,
with the indeterminates $t$ and $u$ being replaced by $1$.
Alternatively, it may be deduced from Proposition~\ref{invhomlocenv}
using standard skew polynomial ring results and
Corollary~\ref{leacor}.
\end{proof}

\begin{rmk}
It is a routine matter to check that, with $A$ as in \ref{invexactP}
and $f=x_1(4-x_3^2)+x_2^2$ as in Example~\ref{Plocenv}, $T$ and $T/gT$ are quantizations of $A$ and
$A/fA$ respectively, that $T_1$ and $T_1/pT_1$ are deformations of $A$ and
$A/fA$ respectively and that the situation is as illustrated in
Figure~\ref{quantPoiss}.
\end{rmk}
\begin{rmk}
\label{altfilt} In the proof of Proposition~\ref{invhomlocenv}, use
was made of a filtration of $T$ for which $\gr (T)$ is commutative.
There are other filtrations for which $\gr (T)$ is a non-commutative
iterated skew polynomial ring and a quantization of another exact
Poisson bracket. If we take the degree function $d$ on $M_4$ such
that
\begin{equation*}
d(z_1)=2,\; d(z_2)=3,\; d(z_3)=2 \mbox{ and } d(u)=1,
\end{equation*}
then \eqref{tcent}, \eqref{X12rel} and \eqref{X13rel} become
homogeneous while, in \eqref{X23rel}, only the term $4t$ has degree
less than $5$. For the $d$-standard filtration of $T$, $\gr (T)$
has, like $T$, the form
$C=\F[t,x_3][x_2;\delta][x_1;\sigma,\delta_1]$ but with
$\delta=-tx_3^2\partial/\partial x_3$ and with the central element
$g=-x_3^2x_1+x_2^2+3tx_3x_2+t^2x_3^2$. Methods similar to those used
in the study of $T$ confirm the existence of such a skew polynomial
ring $C$ and that $\gr (T)\simeq C$. Setting $f=x_2^2-x_3^2x_1\in
A=\F[x_1,x_2,x_3]$, the algebras $\gr (T)$ and $\gr (T/gT)$ are
respectively quantizations of $A$, for the exact Poisson bracket
$\{-,-\}_f$ and the coordinate ring $A/fA$ of the Whitney umbrella.
Alternatively, if we take $d(z_1)=3, d(z_2)=4, d(z_3)=2$ and
$d(u)=1$ and set $h=x_2^2$ then $\gr (T)$ is a quantization of $A$
with the exact Poisson bracket $\{-,-\}_h$ and is isomorphic to the
enveloping algebra of the three-dimensional Heisenberg Lie algebra.
\end{rmk}

\begin{rmk}\label{potential}
The $\F$-algebra $T$ in Proposition~\ref{thinvlocenv} is an example
of an algebra determined by a noncommutative potential. In the
notation of \cite{Ginzburg}, it is $\mathfrak{U}(F,\Phi)$ where $F=\F\langle x_1,x_2,x_3\rangle$ and
\[\Phi=x_1x_2x_3-x_3x_2x_1+x_1x_3^2-x_2x_3-x_2^2-\tfrac{1}{2}x_3^2-4x_1,\] so that
\begin{align*}
\frac{\partial \Phi}{\partial x_2}&=x_3x_1-x_1x_3-x_3-2x_2, \\
\frac{\partial \Phi}{\partial x_1}&=x_2x_3-x_3x_2+x_3^2-4\mbox{ and }\\
\frac{\partial \Phi}{\partial
x_3}&=x_1x_2-x_2x_1+x_3x_1+x_1x_3-x_2-x_3.
\end{align*}
\end{rmk}

\section{Invariants for the quantum torus}
\label{invqtorus} The aim of this section is to complete
Figure~\ref{quantPoiss} for the quantum torus $S_q=W_q$ and the
generic quantum torus $S=W_Q$ and their reversing automorphisms,
$\theta_q$ and $\theta$, specified in Example~\ref{qtorus}. We shall
assume that $q$ is not a root of unity and, in order to apply
standard results on fixed rings, that
 $\ch \F\neq 2=|\langle \theta\rangle|$.

The situation is more complex than in the previous section due to
the invertibility of $Q$ and the lack of any apparent iterated skew
polynomial ring structure for the quantization $T$ or the
deformations, which in this case are parametrized by $q\in
\F\backslash\{0\}$ and will be written $T_q$ rather than $T_{q-1}$.
The ring $T_q$ has been of interest in mathematical physics
\cite{fairlie,mp1,mp2} and in the work of Terwilliger
\cite{terwilliger} and others on Leonard pairs and Askey-Wilson
relations.

\begin{prop}
\label{PBWgenqtorus}  Let $T$ be the $\F$-algebra
generated by $Q$, $Q^{-1}$, $x_1, x_2$ and $x_3$ subject to the
relations
\begin{align}\label{Qcent} x_iQ&=Qx_i,\; x_iQ^{-1}=
Q^{-1}x_i\mbox{ for }i=1,2,3,\; QQ^{-1}=Q^{-1}Q=1,\\
x_1x_2&=Qx_2x_1+(1-Q^2)x_3, \label{QX12}\\
x_2x_3&=Qx_3x_2+(Q^{-1}-Q)x_1,\label{QX23}\\
x_1x_3&=Q^{-1}x_3x_1+(1-Q^{-2})x_2. \label{QX13}
\end{align}
The algebra $T$ has a partially localized PBW basis
$\{Q^ix_3^jx_2^kx_1^\ell: i\in \Z, j,k,\ell\in \N_0\}$.
\end{prop}
\begin{proof}
(i) Let $F$ be the free algebra $F_5$ and $M$ be the free monoid  $M_5$.
It will be convenient to write $u$ for
$z_4$ and $v$ for $z_5$. Let $\psi:F\rightarrow T$ be the
surjective homomorphism such that $\psi(z_i)=x_i$, $1\leq i\leq 3$,
$\psi(u)=Q$ and $\psi(v)=Q^{-1}$. Let $d_1:M_5\rightarrow \N_0$ be the
degree function such that
\[d_1(u)=d_1(v)=d_1(z_1)=0\mbox{ and }d_1(z_2)=d_1(z_3)=1.\]
We shall apply the Diamond Lemma with the augmented
$d_1$-length-lexicographic ordering, as defined in
Definitions~\ref{dfiltrations}. In the $d_1$-standard filtration of
$T$, the largest monomials appearing in \eqref{QX12}, \eqref{QX23}
and \eqref{QX13} are $x_1x_2$, $x_2x_3$ and $x_1x_3$ respectively.
Overlap ambiguities involving the relations in \eqref{Qcent} are
easily resolved. The only other
 ambiguity is the overlap ambiguity $(x_1x_2)x_3=x_1(x_2x_3)$. Reducing
$(x_1x_2)x_3$ by applying \eqref{QX12} followed by \eqref{QX13} and
\eqref{QX23}, together with several applications of \eqref{Qcent},
one obtains
\[Qx_3x_2x_1+(Q^{-1}-Q)x_1^2+(Q-Q^{-1})x_2^2+(1-Q^2)x_3^2.\] The
same result is obtained by reducing $x_1(x_2x_3)$ using
\eqref{QX23}, \eqref{Qcent}, \eqref{QX13} and \eqref{QX12}. It
follows, by the Diamond Lemma, that $T$ has the stated
basis.
\end{proof}

Although the degree function $d_1$ used above will be helpful in
showing that $T/gT$ is a domain, for a central element $g$ to be
specified in Proposition~\ref{gcentral}, we shall also make use of
the degree function $d_2:M_5\rightarrow \N_0$ for which
\[d_2(u)=d_2(v)=0\mbox{ and }d_2(z_1)=d_2(z_2)=d_2(z_3)=1.\]
This has the advantage that, after passage, via localization at
$\F[Q^{\pm1}]\backslash\{0\}$, to a filtered algebra over $\F(Q)$,
the $d_2$-standard-filtration becomes finite. The following Lemma
identifies the associated graded rings for the
$d_i$-standard-filtrations, $i=1,2$.
\begin{lemma}
\label{filtQ} (i) There exist $\sigma\in\Aut_\F(\F[Q^{\pm 1},y_3])$,
$\sigma_1\in\Aut_\F(\F[Q^{\pm 1},y_3][y_2;\sigma])$ and a
$\sigma_1$-derivation $\delta$ of $\F[Q^{\pm 1},y_3][y_2;\sigma]$
such that, for the $d_1$-standard filtration of $T$,
 $\gr(T)$ is an iterated skew polynomial ring $\F[Q^{\pm
1},y_3][y_2;\sigma][y_1;\sigma_1,\delta]$. The algebra $\gr(T)$ is
generated by $Q^{\pm1}, y_1, y_2$ and $y_3$ subject to the
relations:
\begin{gather*}
y_1Q=Qy_1,\quad y_2Q=Qy_2,\quad y_3Q=Qy_3,\quad QQ^{-1}=1=Q^{-1}Q,\\
 y_1y_2=Qy_2y_1+(1-Q^2)y_3,\quad
y_2y_3=Qy_3y_2,\quad y_1y_3=Q^{-1}y_3y_1+(1-Q^{-2})y_2.
\end{gather*}

(ii) For the $d_2$-standard filtration of $T$,
 $\gr(T)$ is an iterated skew polynomial ring
 $\F[Q^{\pm1},y_3][y_2;\sigma][y_1;\sigma_1]$, where $\sigma$  and $\sigma_1$ are as in (i). The
algebra $\gr(T)$ is generated by $Q^{\pm1}, y_1, y_2$ and $y_3$
subject to the relations:
\begin{gather*}
y_1Q=Qy_1,\quad y_2Q=Qy_2,\quad y_3Q=Qy_3,\quad QQ^{-1}=1=Q^{-1}Q,\\
y_1y_2=Qy_2y_1,\quad y_2y_3=Qy_3y_2,\quad
y_1y_3=Q^{-1}y_3y_1.\end{gather*}
\end{lemma}
\begin{proof}
(i) Let $\sigma\in \Aut_\F(\F[Q^{\pm 1},y_3])$  be such that
$\sigma(y_3)=Qy_3$ and $\sigma(Q)=Q$. A method similar to that used
in the proof of Proposition~\ref{invhomlocenv} shows that there
exist $\sigma_1\in \Aut_\F(\F[Q^{\pm 1},y_3][y_2;\sigma])$ and a
$\sigma_1$-derivation $\delta$ of $\F[Q^{\pm 1},y_3][y_2;\sigma]$
such that $\sigma_1(Q)=Q$, $\sigma_1(y_2)=Qy_2$,
$\sigma_1(y_3)=Q^{-1}y_3$, $\delta(Q)=0$, $\delta(y_2)=(1-Q^2)y_3$
and $\delta(y_3)=(1-Q^{-2})y_2$. Let
$B=\F[Q^{\pm1},y_3][y_2;\sigma][y_1,\sigma_1,\delta]$. Then $B$ is
generated by $Q^{\pm1}, y_1, y_2$ and $y_3$ subject to the stated
relations and has basis $\{Q^iy_3^jy_2^ky_1^\ell: i\in \Z,
j,k,\ell\in \N_0\}$.

The degree function $d_1$ is compatible with the reduction scheme
represented by the defining relations for $T$, as shown in
Proposition~\ref{PBWgenqtorus}. The defining relations for $B$ are
satisfied  in $\gr(T)$, with each $y_i$ replaced by $\ov{x_i}$.  It
follows from \cite[Proposition 1]{fddaj} that there is a surjection
$\phi:B\rightarrow \gr(T)$ such that $\phi(y_i)=\ov{x_i}$,
$i=1,2,3$, and $\phi(Q)=Q$. (As $Q\in T_0$, we write $Q$ and
$Q^{-1}$ rather than $\ov Q$ and $\ov {Q^{-1}}$ in $\gr(T)$.) By
Proposition~\ref{PBWgenqtorus} and Proposition~\ref{grbasis}, with
the augmented $d_1$-length-lexicographic ordering, $\gr(T)$ has
basis $\{Q^i\ov{x_3}^j\ov{x_2}^k\ov{x_1}^\ell: i\in \Z, j,k,\ell\in
\N_0\}$ so $\phi$ is an isomorphism.

(ii) The proof is similar to that of (i), but simpler, with the
$\sigma_1$-derivation $\delta$ replaced by $0$.
\end{proof}

\begin{cor}
\label{Tdomain} The algebra $T$ is a domain.
\end{cor}
\begin{proof}
This is immediate from \cite[Proposition 1.6.6(i)]{McCR} and either
part of Lemma~\ref{filtQ}.
\end{proof}

\begin{prop}\label{gcentral}
In $T$, let
$g=x_3x_2x_1-Qx_3^2-Q^{-2}x_2^2-x_1^2+2(1+Q^{-2})$. Then $g$ is a
central element of $T$ and $T/(g-\kappa)T$ is a domain for all $\kappa\in \F$.
\end{prop}

\begin{proof}  Using the defining relations, it can be checked that
\begin{align}
\label{one} x_1x_2^2&=Q^2x_2^2x_1+(1-Q^4)x_3x_2+(Q^2-1)^2x_1,\\
\label{two}x_1x_3^2&=Q^{-2}x_3^2x_1+(Q-Q^{-3})x_3x_2-Q^{-1}(Q-Q^{-1})^2x_1,\\
\label{three}x_2^2x_3&=Q^2x_3x_2^2+(Q^{-1}-Q^3)x_2x_1+(Q^2-1)^2x_3,\\
\label{four}x_1^2x_3&=Q^{-2}x_3x_1^2+(Q-Q^{-3})x_2x_1-(Q-Q^{-1})^2x_3,\\
\label{five}x_2x_1x_3&=x_3x_2x_1+(Q^{-2}-1)x_1^2+(1-Q^{-2})x_2^2\mbox{ and}\\
\label{six}x_1x_3x_2&=x_3x_2x_1+(Q^{-1}-Q)x_3^2+(1-Q^{-2})x_2^2.
\end{align}
It follows routinely that $x_1g=gx_1$ and $x_3g=gx_3$. By
\eqref{QX13}, $g$ commutes with $(1-Q^{-2})x_2$ so, as $Q$ is
central and $T$ is a domain, $g$ is central.

 Let $h=g-\kappa$. We filter $T/hT$ using the $d_1$-standard filtration, where $d_1$ is as
 in the proof of Proposition~\ref{PBWgenqtorus}, and, for
 application of the Diamond Lemma, we use the augmented
$d_1$-length-lexicographic ordering.
 We write $u_i$ for $x_i+hT$, $1\leq i \leq 3$, and,
 with a mild abuse of notation, $Q$ for $Q+hT$.
 Thus the defining relations for
 $T/hT$, each written with the largest term isolated on the left hand side, are:
\begin{align}
\label{UQcent} u_iQ&=Qu_i,\; u_iQ^{-1}=Q^{-1}u_i,\;i=1,2,3,\; QQ^{-1}=1=Q^{-1}Q\\
u_1u_2&=Qu_2u_1+(1-Q^2)u_3, \label{QU12}\\
u_2u_3&=Qu_3u_2+(Q^{-1}-Q)u_1,\label{QU23}\\
u_1u_3&=Q^{-1}u_3u_1+(1-Q^{-2})u_2, \label{QU13}\\
\label{U2squared}
u_2^2&=Q^2u_3u_2u_1-Q^3u_3^2-Q^2u_1^2+2(Q^{2}+1)-\kappa
Q^2.\end{align}  There are no inclusion ambiguities  and the only
overlap ambiguities, apart from those that involve \eqref{UQcent},
are \[(u_1u_2)u_3=u_1(u_2u_3),\; (u_1u_2)u_2=u_1(u_2^2)
\mbox{and}\;u_2(u_2u_3)=(u_2^2)u_3.\] The
 first is resolved as in the proof of
Proposition~\ref{PBWgenqtorus} but two formulae obtained during this
calculation, namely
\begin{multline}
\label{fiveU}u_2u_1u_3=Q^2u_3u_2u_1+(Q^{-2}-Q^2)u_1^2+(Q-Q^3)u_3^2\\+2(Q^2-Q^{-2})+
\kappa(1-Q^2)
 \end{multline} and
  \begin{multline}
\label{sixU}u_1u_3u_2=Q^2u_3u_2u_1+(Q^{-1}-Q^3)u_3^2+(1-Q^{2})u_1^2\\+2(Q^2-Q^{-2})+
\kappa(1-Q^2)
 \end{multline}
are used in resolving the other two ambiguities. We reduce
$(u_1u_2)u_2$ and $u_1(u_2^2)$, beginning the former by applying
\eqref{QU12} and the latter by applying \eqref{U2squared}, and
making use of \eqref{QU12}, \eqref{QU23}, \eqref{QU13},
\eqref{fiveU} and \eqref{sixU}, to obtain the result
\[
Q^4u_3u_2u_1^2-Q^5u_3^2u_1-Q^4u_1^3+(1-Q^4)u_3u_2+((3-\kappa)Q^4+1)u_1
\] in both cases.
Similarly, $u_2(u_2u_3)$ and $(u_2^2)u_3$ both reduce to
\[
Q^4u_3^2u_2u_1-Q^4u_3u_1^2-Q^5u_3^3+(Q^{-1}-Q^3)u_2u_1+((3-\kappa)Q^4+1)u_3.
\]
By the Diamond Lemma, $\{Q^iu_3^ju_2^ku_1^\ell:i\in\Z, j,k,\ell\geq
0, k<2\}$ is a basis for $T/hT$.

We claim that $\gr(T/hT)\simeq \gr(T)/\ov{h}\gr (T)$. For $i=1,2,3$,
we write $v_i=\ov {u_i}\in \gr(T/hT)$ (and we write $Q$ for $\ov
Q$). By Proposition~\ref{grbasis}, $\{Q^iv_3^jv_2^kv_1^\ell:i\in\Z,
j,k,\ell\geq 0, k<2\}$ is a basis for $\gr(T/hT)$ whose generators
$Q^{\pm 1}, v_1, v_2$ and $v_3$ satisfy the following relations:
\begin{align*}
v_iQ&=Qv_i,\; v_iQ^{-1}=Q^{-1}v_i,\;i=1,2,3,\; QQ^{-1}=1=Q^{-1}Q,\\
v_1v_2&=Qv_2v_1+(1-Q^2)v_3, \\
v_2v_3&=Qv_3v_2,\\
v_1v_3&=Q^{-1}v_3v_1+(1-Q^{-2})v_2,\\
v_2^2&=Q^2v_3v_2v_1-Q^3v_3^2.
\end{align*}
The same relations, with each $v_i$ replaced by $w_i:=y_i+\ov{h}\gr
(T)\in \gr(T)/\ov{h}\gr(T)$, are defining relations for
$\gr(T)/\ov{h}\gr(T)$ so there is a surjection $\phi:\gr
(T)/\ov{h}\gr(T)\rightarrow \gr(T/hT)$ such that $\phi(Q)=Q$ and
$\phi(w_i)=v_i$, $i=1,2,3$.
 Effectively the same Diamond Lemma calculations as for $T/hT$, but with some
low degree terms deleted, show that $\gr(T)/\ov{h}\gr(T)$ has basis
$\{Q^iw_3^jw_2^kw_1^\ell:i\in\Z, j,k,\ell\geq 0, k<2\}$. It follows
that $\phi$ is an isomorphism.

By \cite[Proposition 1.6.6(i)]{McCR}, it now suffices to show that
$\gr(T)/\ov{h}\gr(T)$ is a domain. For this, recall, from
Lemma~\ref{filtQ}(i), that $\gr(T)=\F[Q^{\pm
1},y_3][y_2;\sigma][y_1;\sigma_1,\delta]$. We apply
\cite[Proposition 1]{normal} to the central element
$\ov{h}=y_3y_2y_1-Qy_3^2-Q^{-2}y_2^2$ of degree $1$ in $y_1$. Let
$D:=\F[Q^{\pm 1},y_3][y_2;\sigma]$, $d=y_3y_2$ and
$e=-Qy_3^2-Q^{-2}y_2^2$. Here $\sigma(y_3)=Qy_3$ and $\sigma(Q)=Q$.
To conclude that $\gr(T)/\ov{h}\gr(T)$ is a domain, we need to check
that $e$ is regular modulo the ideal $Dd$. In $D$, $Dd$ is the
intersection of two height one primes $Dy_3$ and $Dy_2$ which are
completely prime. Hence all zero-divisors modulo $Dd$ are in
$Dy_3\cup Dy_2$ and so $e$ is regular modulo $Dd$. This completes
the proof that $T/hT$ is a domain.
\end{proof}

\begin{prop}
\label{invgenqtorus}Let $S=W_Q$ and the reversing
automorphism $\theta$ be as in Example~\ref{qtorus}(ii). Then $S^\theta$ is isomorphic to $T/gT$.
\end{prop}
\begin{proof}
In the notation of Lemma~\ref{srelations}, let
$a_1=s_0(y)=s_0(y^{-1})=y+y^{-1}$, let $a_2=s_1(1)=x+x^{-1}$ and let
$a_3=s_1(y)=yx+y^{-1}x^{-1}$. Note that
$s_0(\alpha^{-1}(y^{-1}))=s_0(Qy^{-1})=Qa_1$ and
$\gamma(y^{-1})-\alpha^2\gamma(y^{-1})=(1-Q^2)y$. By \eqref{sr} with
$r=y^{-1}$,
\[
a_1a_2-Qa_2a_1=(1-Q^2)a_3.
\]
By \eqref{sr2}, $a_3a_2-Q^{-1}a_2a_3=(1-Q^{-2})a_1$ so
\[
a_2a_3-Qa_3a_2=(Q^{-1}-Q)a_1.
\]
Applying \eqref{rrprime} with $r=r^\prime=y$, we obtain
$a_1a_3-Q^{-1}a_3a_1=s_1((y^{-1}-Q^{-2}y^{-1})y)=(1-Q^{-2})a_2$,
whence
\[
a_3a_1-Qa_1a_3=(Q^{-1}-Q)a_2.
\]

By Corollary~\ref{qtcor}, $a_1, a_2, a_3, Q$ and $Q^{-1}$ generate
$S^{\theta}$ so there is a surjective $\F$-homomorphism
$\eta:T\rightarrow S^\theta$, such that
$\eta(x_i)=a_i$, $i=1,2,3$, and $\eta(Q)=Q$. Therefore $S^\theta\simeq T/\ker \eta$.

Note that $a_1^2=y^2+2+y^{-2}$, $a_2^2=x^2+2+x^{-2}$,
$a_3^2=Qy^2x^2+2Q^{-1}+Qy^{-2}x^{-2}$ and
\begin{align*}
a_3a_2a_1&=Q^2y^2x^2+Q^2y^{-2}x^{-2}+Q^{-2}x^2+Q^{-2}x^{-2}+y^2+y^{-2}+2\\
&=Qa_3^2+Q^{-2}a_2^2+a_1^2-2(1+Q^{-2}).
\end{align*}
Thus $g\in \ker\eta$.

As in the proof of Proposition~\ref{invhomlocenv}, we use
GK-dimension to show that $\ker \eta=gT$. However it will
be convenient to work over the rational function field $\K:=\F(Q)$
rather than over $\F$. To this end, let $\mathcal{C}$ denote the
central multiplicatively closed set $\F[Q^{\pm1}]\backslash\{0\}$.
It follows from Proposition~\ref{PBWgenqtorus} that the localization
$\hat{T}$ of $T$ at $\mathcal{C}$ is a $\K$-algebra with a
PBW basis $\{x_3^ix_2^jx_1^k\}$. As $d_2(Q)=0$, the filtration in
Lemma~\ref{filtQ}(ii) extends to a filtration of $\hat{T}$ as a
$\K$-algebra and, as each $d_2(x_i)>0$, this filtration is finite. The
associated graded $\K$-algebra $\gr (\hat{T})$ is the quantum
coordinate ring of $\K^3$ generated by $y_1,y_2,y_3$ subject to the
relations
\[
y_1y_2=Qy_2y_1,\quad y_2y_3=Qy_3y_2,\quad y_3y_1=Qy_1y_3.
\] If $V$ is the generating subspace $\K y_1+\K y_2+\K y_3$ then $\dim_\K V^n$ is the same as for the commutative polynomial ring $\K[y_1,y_2,y_3]$ so
$\GK_\K \gr(\hat{T})=3$ and, by \cite[Proposition 6.6]{KL} or
\cite[Proposition 8.6.5]{McCR}, $\GK_\K (\hat{T})=3$. Let $\hat{S}$
denote the localization of $S$ at $\mathcal{C}$, so that,
as $\theta(\mathcal{C})=\mathcal{C}$, $\theta$ extends  to $\hat{S}$
 in the obvious way and the surjective
homomorphism $\eta:T\rightarrow S^{\theta}$
extends to a surjective homomorphism $\hat{\eta}:\hat{T}\rightarrow
\hat{S}^{\theta}$.

Now  $\GK_\K(\hat{S}^\theta)=\GK(\hat{S})=2$, by \cite[Corollary
26.13(ii)]{passman} and \cite[Proposition 8.2.9]{McCR}. As
$\hat{S}^{\theta}$ is a domain, $\ker \hat{\eta}$ is a prime ideal
$P$, say, of $\hat{T}$ and, by \cite[Corollary 3.16]{KL},
\[3=\GK_\K(\hat{T})\geq\GK_\K(\hat{T}/P)+\hgt(P)=\GK_\K(\hat{S}^\theta)+\hgt(P)=2+\hgt(P).\]
Hence $\hgt(P)\leq 1$. As $P\neq 0$ and $\hat{T}$ is a domain,
$\hgt(P)=1$. Also $\ker \eta=P\cap T$ so $\ker \eta$ is a prime
ideal of $T$ of height one by \cite[Proposition 2.1.16(vii)]{McCR}.
By Proposition~\ref{gcentral}, $T/gT$ is a domain so, as $g\in
\ker\eta$, it follows that $\ker \eta= gT$ and hence that
$T/gT\simeq S^{\theta}$.
\end{proof}

The following result, identifying the invariants for the reversing
automorphism of the quantum torus in Example~\ref{qtorus} rather
than the generic quantum torus, can be proved either by adapting the
methods above, with simplification due to the replacement of the
invertible indeterminate $Q$ by the non-zero scalar $q$, or by
applying the results above together with Corollary~\ref{qtcor}. The
algebra $T_q$ defined in the statement is isomorphic to $T/(Q-q)T$
and $p$ is the image of $g$ in $T_q$.
\begin{prop}
\label{thinvqtorus} Let $S_q=W_q$ and the reversing
automorphism $\theta_q$ be as in Example~\ref{qtorus}(i).
Let $T_q$ be the $\F$-algebra generated by $ x_1, x_2$ and $x_3$
subject to the relations
\begin{align}
x_1x_2-qx_2x_1&=(1-q^2)x_3,\label{qx12} \\
x_2x_3-qx_3x_2&=(q^{-1}-q)x_1,\label{qx23} \\
x_3x_1-qx_1x_3&=(q^{-1}-q)x_2.\label{qx31}  \end{align}

(i) $T_q$ has a PBW basis $\{x_3^ix_2^jx_1^k:i,j,k\geq 0\}$.

(ii) $T_q$ has a filtration for which $\deg x_1=0$ and $\deg
x_2=\deg x_3=1$ and the associated graded ring is an iterated skew
polynomial ring over $\F$ generated by $y_1, y_2$ and $y_3$ subject
to the relations
\begin{equation*}
 y_1y_2-qy_2y_1=(1-q^2)y_3, \quad
y_2y_3-qy_3y_2=0,\quad y_3y_1-qy_1y_3=(q^{-1}-q)y_2. \end{equation*}

(iii) $T_q$ has a filtration for which $\deg x_1=\deg x_2=\deg x_3=1$ and the associated graded ring
is an iteration skew polynomial ring over $\F$ generated by
$y_1, y_2$ and $y_3$ subject to the relations
\begin{equation*}
y_1y_2=qy_2y_1, \quad y_2y_3=qy_3y_2,\quad y_3y_1=qy_1y_3.
\end{equation*}

(iv) Let $p=x_3x_2x_1-qx_3^2-q^{-2}x_2^2-x_1^2+2(1+q^{-2})$. Then
$p$ is a central element of $T_q$ and $T_q/pT_q$ is a domain.

(v) $S_q^{\theta_q}$ is isomorphic to $T_q/pT_q$.
\end{prop}

\begin{rmk}
Setting $t=Q-1$, which is a central regular non-unit, $T/tT\simeq \F[x_1,x_2,x_3]$, the commutative polynomial
algebra. It is a routine matter to confirm that, in accordance with
the discussion in Example~\ref{qtorus}(ii), $T$ is a
quantization of the Poisson algebra $A=\F[x_1,x_2,x_3]$ with
\begin{equation*}
\{x_1,x_2\}=x_1x_2-2x_3, \quad \{x_2,x_3\}=x_2x_3-2x_1,\quad
\{x_3,x_1\}=x_1x_3-2x_2,
\end{equation*}
that is with the exact Poisson bracket determined by
$f=x_1x_2x_3-x_1^2-x_2^2-x_3^2+4$. Also $T/gT$ is a quantization of
the Poisson algebra $A/fA$.
\end{rmk}
\begin{rmk}
The associated graded rings in Lemma~\ref{filtQ} are also
quantizations of $A$ with appropriate exact Poisson brackets
$\{-,-\}_f$. In Lemma~\ref{filtQ}(i) $f=x_1x_2x_3-x_2^2-x_3^2$ while
in Lemma~\ref{filtQ}(ii), $f=x_1x_2x_3$. In both cases $\gr(T/gT)$
is a quantization of the Poisson algebra $A/fA$, with
$\gr(T_q/pT_q)$ as a deformation. Note the cyclic form of the
relations in Lemma~\ref{filtQ}(ii) and
Proposition~\ref{thinvqtorus}(iii). In the latter, this
distinguishes $\gr (T_q)$ from quantum affine space $\mathcal{O}_q$
as defined in \cite[p. 15]{BGl}.

Another filtration giving rise to an associated graded ring that is
a quantization of $A$ with an exact Poisson bracket is obtained by
taking $Q$ and $Q^{-1}$ to have degree $0$, $x_2$ and $x_3$ to have
degree $1$ but $x_1$ to have degree $2$. For this filtration, $\gr(T)$ is generated by $Q^{\pm 1}$ and $y_i$, $i=1,2,3$,
subject to the relations
\begin{gather*}y_iQ=Qy_i,\; y_iQ^{-1}=Q^{-1}y_i\mbox{ for }i=1,2,3,\; QQ^{-1}=Q^{-1}Q=1,\\
y_1y_2=Qy_2y_1, \quad y_2y_3=Qy_3y_2+(Q^{-1}-Q)y_1,\quad
y_1y_3=Q^{-1}y_3y_1.
\end{gather*}
There is a corresponding filtration of $T_q$ with $\gr (T_q)\simeq
\gr(T)/(Q-q)\gr(T)$.  Here $\gr(T)$ and $\gr(T_q)$ are respectively
a quantization and deformation of $A$ for the exact Poisson bracket
$\{-,-\}_f$ where $f=x_1x_2x_3-x_1^2$. For each $q\in \F$, $\gr
{T_q}$ is an iterated skew polynomial ring over $\F$ and is an
ambiskew polynomial ring in the sense of \cite{ambi} and, by the
results of \cite{ambi}, it is a down-up algebra in the sense of
\cite{downup}.

Note that, for two of these filtrations, the polynomial $f$ is
reducible and $\gr(T/gT)$ and $\gr(T_q/pT_q)$
are not domains.
\end{rmk}

\begin{rmk}\label{isos}
There is some flexibility in the defining relations for $T_q$. Let
$a,b,c\in\F\backslash\{0\}$ and let $T_q(a,b,c)$ denote the
$\F$-algebra generated by $x_1, x_2$ and $x_3$ subject to the
relations
\begin{equation*}
x_1x_2-qx_2x_1=ax_3,\quad x_2x_3-qx_3x_2=bx_1,\quad
x_3x_1-qx_1x_3=cx_2.
\end{equation*}
Thus $T_q=T_q(1-q^2,q^{-1}-q,q^{-1}-q)$. It can be checked that
there exist $\lambda_i\in \F\backslash\{0\}$, $i=1,2,3$, such that
there is an isomorphism $\theta:T_q(1,1,1)\rightarrow T_q(a,b,c)$
with $\theta(x_i)=\lambda_i x_i$. Algebras isomorphic to
$T_q(1,1,1)$ have been considered, sometimes with $q^2$ in place of
$q$, in the mathematical physics literature \cite{fairlie,mp1,mp2}
and in the literature on Leonard triples and Askey-Wilson relations,
for example \cite{terwilliger,terwilliger+}. The Askey-Wilson
relations for $T_q$ are
\begin{align*}
(1+q^2)x_2x_1x_2-qx_2^2x_1-qx_1x_2^2+q^2x_2x_1x_2&=q^{-1}(1-q^2)^2x_1\mbox{ and }\\
(1+q^2)x_1x_2x_1-qx_1^2x_2-qx_2x_1^2+q^2x_1x_2x_1&=q^{-1}(1-q^2)^2x_2.
\end{align*}
and are obtained by using \eqref{qx12} to substitute for $x_3$, in
terms of $x_1$ and $x_2$, in \eqref{qx23} and \eqref{qx31}.

Let us extend our definitions of quantization and deformations to
noncommutative algebras by defining a quantization of a
noncommutative algebra $A$ to be an algebra $B$ with a central,
regular nonunit $t$ such that $B/tB\simeq A$ and a deformation of
$A$ to be an algebra of the form $B/(t-\lambda)B$ for some
quantization $B$ and some $\lambda\in \F$ for which $t-\lambda$ is a
nonunit. Then $T_q(1,1,1)$ is a deformation of the enveloping
algebra $U$ of the Lie algebra $so_3$, the appropriate
quantization being the $\F$-algebra generated by $Q^{\pm 1}$, $x_1,
x_2$ and $x_3$ subject to the relations
\begin{gather*} x_iQ=Qx_i,\; x_iQ^{-1}=
Q^{-1}x_i\mbox{ for }i=1,2,3,\; QQ^{-1}=Q^{-1}Q=1,\\
x_1x_2-Qx_2x_1=x_3, \quad x_2x_3-Qx_3x_2=x_1,\quad
x_3x_1-Qx_1x_3=x_2.
\end{gather*}

Note that if $q^2\neq 1$, $T_q=T_q(1-q^2, q^{-1}-q, q^{-1}-q)$ is a
deformation of the commutative polynomial algebra $\F[x_1,x_2,x_3]$
and, being isomorphic to
 $T_q(1,1,1)$, it is also a deformation of $U$. This gives rise to a
 dichotomy in its behaviour. For example,
the trivial representation, with $x_1,x_2,x_3$ each acting as $0$,
can be viewed as a deformation of the unique one-dimensional
representation of $U$. On the other hand, if $G$ is the Klein
$4$-group $\{e,a_1,a_2,a_3\}$, there is a surjection
$\phi:T_q(1-q,1-q,1-q)\rightarrow\F G$ with $\phi(x_i)=a_i$,
$i=1,2,3.$ It follows that $T_q$ has a further four $1$-dimensional
representations. This typifies the finite-dimensional simple
representations of $T_q$. For each $d\geq 1$, Fairlie \cite{fairlie}
constructed a $d$-dimensional simple representation of $T_q$ while
Havlicek, Klimyk and Posta constructed another four in \cite{mp1}
and, in \cite{mp2}, Havlicek and Posta showed that there were no
more. The finite-dimensional simple representations of $T_q$ are
also classified, by an independent method, in \cite{nongthesis}.
\end{rmk}

\begin{rmk}
When $T_q$ is viewed as a deformation of the Poisson algebra
$A=\F[x_1,x_2,x_3]$ in Example~\ref{qtorus}(ii), the central element
$p$ is a deformation of the element $x_1x_2x_3-x_1^2-x_2^2-x_3^2+4$
of $\PZ(A)$. On the other hand, when $T_q$ is viewed as a
deformation of the enveloping algebra $U$ of the Lie algebra $so_3$,
the cubic term of $p$ has a coefficient of the form $g(q)$ where
$g(Q)\in \F[Q^{\pm 1}]$ and $g(1)=0$ so, although  $p$ is cubic, it
is a deformation of a quadratic element of $U$ which, up to scalar
multiplication and translation, is the Casimir element of $U$.
\end{rmk}

\begin{rmk}
In common with $T$ in Remark~\ref{potential}, the $\F$-algebra $T_q$
in Proposition~\ref{thinvlocenv} is  determined by a noncommutative
potential. Here it is a variation on one of the basic examples of
such an algebra, see \cite[Example 1.3.6]{Ginzburg},
$T_q=\mathfrak{U}(\F\langle x_1,x_2,x_3\rangle,\Pi_q)$ where
\[\Pi_q=x_1x_2x_3-qx_3x_2x_1+\tfrac{1}{2}(q-q^{-1})(x^2+y^2+qz^2).\]
\end{rmk}

\begin{rmk}
 By Proposition~\ref{thinvqtorus}(v), the
quotient division ring $Q(T_q/pT_q)$ is isomorphic to
$Q(W_q^{\theta_q})$ which, by \cite[Theorem 10.5.19(v)]{McCR}, is equal
to the ring of invariants $Q(W_q)^{\theta_q}$ for the induced action of
$\theta_q$ on the quotient division ring $Q(W_q)$ of the quantum
plane. In \cite[13.6]{jtsmvdb}, Stafford and Van den Bergh have
shown that, if $q$ is not a root of unity, $Q(W_q)^\theta\simeq
Q(W_q)$ and have presented, without details of the calculation, a
pair of generators $f,g$ for $Q(W_q)^{\theta_q}$, as a division algebra,
satisfying $fg=qgf$. The elements of $Q(T_q/pT_q)$ corresponding to
$f$ and $g$ are $(a_3a_2-2a_1)(a_2^2-4)^{-1}$ and
$(2a_3-a_1a_2)(a_2^2-4)^{-1}$ respectively, where each $a_i$ is the
image of $x_i$ in $T_q/pT_q$.
\end{rmk}

\end{document}